\documentclass[final]{siamltex}
\usepackage{subfigure}

\usepackage[usenames]{color}
\usepackage[color]{showkeys}
\definecolor{refkey}{gray}{0.8}
\definecolor{labelkey}{gray}{0.8}

\usepackage[left=1in,top=1in,right=1in,bottom=1in,letterpaper]{geometry}
\usepackage{setspace,url}
\usepackage{verbatim}
\usepackage{longtable}

\usepackage[leqno]{amsmath}
\usepackage{amsxtra, amsfonts,amscd, amssymb, graphicx,epsf}
\usepackage{algorithm,algorithmic}
\usepackage[algo2e, vlined,ruled]{algorithm2e}

\numberwithin{equation}{section}

\newcommand{\be}{\begin{equation}}
\newcommand{\ee}{\end{equation}}
\newcommand{\ba}{\begin{array}}
\newcommand{\ea}{\end{array}}
\newcommand{\bea}{\begin{eqnarray}}
\newcommand{\eea}{\end{eqnarray}}
\newcommand{\beaa}{\begin{eqnarray*}}
\newcommand{\eeaa}{\end{eqnarray*}}

\newcommand{\half}{\frac{1}{2}}
\newcommand{\br}{\mathbb{R}}

\usepackage[usenames]{color}

\newtheorem{remark}[theorem]{Remark}

\newcounter{alg-counter} 

\usepackage{multirow}
\newcommand{\bean}{\begin{eqnarray}\nonumber}
\newcommand{\eproof}{$\quad \Box$}
\newcommand{\dis}{\displaystyle}
\newcommand{\df}[2]{\dis{\frac{#1}{#2}}}
\newcommand{\A}{\mathcal{A}}
\newcommand{\X}{\mathcal{X}}
\newcommand{\zero}{\mathbf{0}}
\newcommand{\Diag}{\mbox{Diag}}
\newcommand{\vvec}{\mbox{vec}}
\newcommand{\Tr}{\mbox{Tr}}
\newcommand{\etal}{{ et al. }}

\newcommand{\st}{\mbox{ s.t.}}
\newcommand{\mtn}{{m\times n}}

\newcommand{\myspan}{\mbox{span}}
\begin{document}
\title{Convergence of fixed-point continuation algorithms for \\ matrix rank minimization}
\author{Donald Goldfarb\footnotemark[2] \ and Shiqian Ma\footnotemark[2]}

\renewcommand{\thefootnote}{\fnsymbol{footnote}}
\footnotetext[2]{Department of Industrial Engineering and Operations Research,
Columbia University, New York, NY 10027. Email:
\{goldfarb,sm2756\}@columbia.edu. Research supported in part by
NSF Grants DMS 06-06712 and DMS 10-16571, ONR Grants N00014-03-0514 and
N00014-08-1-1118, and DOE Grants DE-FG01-92ER-25126 and
DE-FG02-08ER-25856.}

\renewcommand{\thefootnote}{\arabic{footnote}}

\maketitle

\begin{center} {June 18, 2009. \quad  This version: December 28, 2010} \end{center}


\begin{abstract} The matrix rank minimization problem has applications in many fields
such as system identification, optimal control, low-dimensional
embedding, etc. As this problem is NP-hard in general, its
convex relaxation, the nuclear norm minimization problem, is often
solved instead. Recently, Ma, Goldfarb and Chen proposed a fixed-point continuation algorithm for solving the nuclear norm
minimization problem \cite{Ma-Goldfarb-Chen-2008}. By incorporating
an approximate singular value decomposition technique in this
algorithm, the solution to the matrix rank minimization problem is
usually obtained. In this paper, we study the
convergence/recoverability properties of the fixed-point
continuation algorithm and its variants for matrix rank
minimization. Heuristics for determining the rank of the matrix when
its true rank is not known are also proposed. Some of these
algorithms are closely related to greedy algorithms in compressed
sensing. Numerical results for these algorithms for solving affinely
constrained matrix rank minimization problems are reported.
\end{abstract}

\begin{keywords} Matrix Rank Minimization, Matrix Completion, Greedy Algorithm, Fixed-Point Method, Restricted Isometry Property, Singular Value Decomposition
\end{keywords}

\begin{AMS} Primary, 90C59; Secondary, 15B52, 15A18 \end{AMS}

\section{Introduction} In this paper, we are interested in the affinely constrained matrix rank minimization (MRM)
problem, which
can be cast as \bea\ba{ll}\label{MRM-affine-rank}\min & \rank(X) \\
\st & \A(X)=b,\ea\eea where $X\in\mathbb{R}^{\mtn}$,
$b\in\mathbb{R}^p$ and $\A:\mathbb{R}^{\mtn}\mapsto\mathbb{R}^p$ is
a linear map. Without loss of generality, we assume that $m\leq n$ throughout this paper.


Problem \eqref{MRM-affine-rank} 
has applications in many fields such as system identification
\cite{Liu-Vandenberghe-2008}, optimal control
\cite{ElGhaoui-Gahinet-1993, Fazel-Hindi-Boyd-2001,
Fazel-Hindi-Boyd-2004}, and low-dimensional embedding in Euclidean
space \cite{Linial-London-Rabinovich-1995}, etc. For example,
consider the problem of designing a low-order discrete-time
controller for a plant, so that the step response of the combined
controller and plant lies within specified bounds. Suppose the plant
impulse response is $h(t),t=0,\ldots,N$, the controller impulse
response is $x(t),t=0,\ldots,N$, and $u(t)=1, t=0,\ldots,N$ is the
step input. Then finding a low-order system is equivalent to solving
the following problem: \bea\label{prob:low-order-system}\ba{ll}\min
& \rank(\mathcal{H}(x))
\\ \st & b_l(t)\leq (h*x*u)(t)\leq b_u(t), t=0,\ldots,N,\ea\eea
where $b_l$ and $b_u$ are given lower and upper bounds on the step
response, $*$ denotes the convolution operator, and $\mathcal{H}(x)$
is the Hankel matrix (see e.g.,
\cite{Fazel-Hindi-Boyd-2003,Sontag-book-1998}). Problem
\eqref{prob:low-order-system} is an application of an inequality-constrained variant of
\eqref{MRM-affine-rank}.

A special case of \eqref{MRM-affine-rank} is the matrix completion problem: \bea\ba{ll}\label{MRM-MatCompletion-rank}\min & \rank(X) \\
\st & X_{ij} = M_{ij}, \quad \forall (i,j)\in\Omega.\ea\eea This problem
has applications in online recommendation systems, collaborative
filtering \cite{Srebro-thesis-2004,Srebro-Jaakkola-2003}, etc., including the
famous Netflix problem \cite{netflixprize}. In the latter problem, users provide ratings to some of
the movies in a list of movies. Here $M_{ij}$ is the rating given to
$j$-th movie by the $i$-th user. Since users only rate a limited number
of movies in the list, we only know some of the entries of the
matrix $M$. The goal of the Netflix problem is to fill in the
missing entries in this matrix. It is commonly believed that only a
few factors contribute to people's tastes in movies. Thus the matrix
$M$ will generally be of low rank. Finding this low-rank completion
to $M$ is just the matrix completion problem
\eqref{MRM-MatCompletion-rank}.

If $X$ is a diagonal matrix, then \eqref{MRM-affine-rank} becomes
the compressed sensing problem \cite{Candes-Romberg-Tao-2006,Donoho-06}: \bea\ba{ll}\label{compressed-sensing-L0} \min &
\|x\|_0 \\ \st & Ax=b, \ea\eea where $A\in\br^\mtn, b\in\br^m$, and
$\|x\|_0$, which is called the $\ell_0$ norm, counts the number of
nonzero elements in the vector $x$. The compressed sensing problem,
which is currently of great interest in signal processing, is
NP-hard \cite{Natarajan-95}. Recent results in compressed sensing
have shown that under certain randomness hypotheses, the optimal solution to
\eqref{compressed-sensing-L0} can be found by solving a convex relaxation of \eqref{compressed-sensing-L0} using only a limited number of
measurements. Since the convex envelope of the function
$\|x\|_0$ on the set $\{x\in\br^n:  \|x\|_\infty\leq 1\}$ is the $\ell_1$ norm
$\|x\|_1:=\sum_{i}|x_i|$ \cite{Hiriart-Urruty-Lemarechal-1993}, a natural choice for a convex relaxation of
problem \eqref{compressed-sensing-L0} is the problem:
\bea\ba{ll}\label{compressed-sensing-L1} \min & \|x\|_1 \\ \st &
Ax=b. \ea\eea Many algorithms for solving
\eqref{compressed-sensing-L0} and \eqref{compressed-sensing-L1} have
been proposed. These include greedy algorithms
\cite{Tibshirani-96,Donoho-Tsaig-Drori-Starck-06,Tropp-06,Donoho-Tsaig-2006,Needell-Tropp-2008,Dai-Milenkovic-08,Blumensath-Davies-08,Blumensath-Davies-2009} for
\eqref{compressed-sensing-L0} and convex optimization algorithms
\cite{Candes-Romberg-2005-l1-magic,Figueiredo-Nowak-Wright-2007,Hale-Yin-Zhang-SIAM-2008,Kim-Koh-Lustig-Boyd-Gorinevsky-2007,vandenBerg-Friedlander-2008,YinOsherGoldfarbDarbon2008}
for \eqref{compressed-sensing-L1}. See
\cite{rice-compressed-sensing-site} for more information on the
theory and algorithms for compressed sensing.

The matrix rank minimization problem
\eqref{MRM-affine-rank} is also NP-hard. To get a tractable problem,
we can replace $\rank(X)$ by the
nuclear norm $\|X\|_*$ of $X$, the convex envelope of $\rank(X)$ on the set $\{X\in\br^\mtn: \|X\|_2\leq 1\}$ \cite{Recht-Fazel-Parrilo-2007}, as proposed by Fazel \etal
\cite{Fazel-Hindi-Boyd-2001}. The nuclear norm of $X$ is defined as the sum of the nonzero singular values of
$X$ and the spectral norm $\|X\|_2$ is equal to the largest singular value of $X$; i.e., if the singular values of $X$ are
$\sigma_1\geq\sigma_2\geq\ldots\geq\sigma_r>\sigma_{r+1}=\ldots=\sigma_{m}=0,$
then
\begin{align*}\|X\|_*=\sum_{i=1}^r\sigma_i\end{align*} and $\|X\|_2=\sigma_1$.
Thus, the nuclear norm relaxation of \eqref{MRM-affine-rank} is:
\bea\ba{ll}\label{MRM-affine-nuclear}\min & \|X\|_* \\ \st &
\A(X)=b. \ea\eea

Let $A$ be the matrix version of $\A$, i.e., $\A(X)=A\cdot \vvec(X)$, where $\vvec(X)$ is the vector obtained by stacking the columns of the matrix $X$ in natural order. Recht \etal \cite{Recht-Fazel-Parrilo-2007} proved that
if the entries of $A$ are drawn from some random distribution and the number
of measurements $p\geq Cr(m+n)\log(mn)$, then with very high
probability, most $m\times n$ matrices of rank $r$ can be recovered
by solving problem
\eqref{MRM-affine-nuclear}, where $C$ is a positive constant; i.e.,
an optimal solution to \eqref{MRM-affine-nuclear} gives an optimal
solution to \eqref{MRM-affine-rank}.

If $b$ is contaminated by noise, then \eqref{MRM-affine-nuclear}
should be relaxed to \bea\ba{ll}\label{MRM-affine-noise-b}\min &
\|X\|_*
\\ \st & \|\A(X)-b\|_2 \leq \theta, \ea\eea where $\theta>0$ is the
noise level. The Lagrangian version of \eqref{MRM-affine-noise-b}
can be written as \bea\label{MRM-affine-unconstrained}\min
\mu\|X\|_*+\half\|\A(X)-b\|_2^2,\eea where $\mu$ is a Lagrangian
multiplier.

Several algorithms have been proposed for solving
\eqref{MRM-affine-rank} and \eqref{MRM-affine-nuclear}. Using the
fact that \eqref{MRM-affine-nuclear} is equivalent to the
semidefinite programming (SDP) problem
\bea\label{MRM-affine-nuclear-SDP-primal}\ba{cl}
\displaystyle\min_{X,W_1,W_2}& \half(\Tr(W_1)+\Tr(W_2)) \\\st &
\begin{bmatrix}W_1 & X \\ X^\top & W_2\end{bmatrix}\succeq 0, \\ & \A(X)=b,
\ea\eea where $\Tr(W)$ denotes the trace of
the square matrix $W$, Recht, Fazel and Parrilo
\cite{Recht-Fazel-Parrilo-2007} and Liu and Vandenberghe
\cite{Liu-Vandenberghe-2008} proposed interior-point methods to
solve this SDP. However, these interior-point methods cannot be used to
solve large problems. First-order methods were proposed by Cai,
Cand\`es and Shen \cite{Cai-Candes-Shen-2008} and Ma, Goldfarb and
Chen \cite{Ma-Goldfarb-Chen-2008} that can solve very large matrix
rank minimization problems efficiently. One of the
algorithms in \cite{Ma-Goldfarb-Chen-2008}, which is called FPCA
(Fixed-Point Continuation with Approximation SVD), almost always
achieves the best recoverability. FPCA can recover $\mtn$ matrices
of rank $r$ using $p$ samples even when $r$ is very close to
the largest rank $r_{\max}:=\max\{r|r(m+n-r)/p<1\}$ of $\mtn$ matrices that one can
recover with only $p$ samples. In this paper, we study the
convergence/recoverability properties and numerical performance of FPCA and some of its
variants. Our main contribution is a weakening of the conditions previously given by Lee and Bresler \cite{Lee-Bresler-2009,Lee-Bresler-ADMIRA-fullpaper-2009} required for the approximate recovery of a low-rank matrix.

{\bf Notation.} We use $\mathbb{R}_+^n$ to denote the nonnegative
orthant of $\mathbb{R}^n$.
We use $\A^*$ to denote the adjoint operator of $\A$. We define the
inner product of two matrices $X$ and $Y\in\mathbb{R}^\mtn$ to be
$\langle X,Y \rangle=\Tr(X^\top Y)=\Tr(Y^\top X)$, and denote the Frobenius norm of the matrix $X$ by
$\|X\|_F=(\Tr(X^\top X))^{1/2}$
and the Euclidean norm of the vector $x$ by $\|x\|_2$. Henceforth, we will write $\A(X)$ as $\A X$ as this should not cause any confusion. For example, $\A^*\A X:= \A^*(\A(X))$.

{\bf Outline.} The rest of this paper is organized as follows. In
Section 2 we review the role that the restricted isometry property
plays in the theory of compressed sensing and matrix rank
minimization. We also present three propositions from
\cite{Lee-Bresler-2009} that provide the basis for the theoretical
results that we give later in the paper. We review the Fixed-Point Continuation (FPC) and FPC with Approximate SVD (FPCA)
algorithms proposed in \cite{Ma-Goldfarb-Chen-2008} in Section 3.
We then address the first variant of FPCA, which we call iterative
hard thresholding (IHT), and prove convergence results for it in
Section 4. Section 5 is devoted to another variant of FPCA, which is
called iterative hard thresholding with matrix shrinkage (IHTMS),
and convergence results for it. We establish
convergence/recoverability properties of FPCAr, a very close variant of FPCA, in Section 6. Some
practical issues regarding numerical difficulties and ways to
overcome them are discussed in Section 7. Finally, we give some
numerical results obtained by applying these algorithms to both randomly
created and real matrix rank minimization problems in Section 8.

%

\section{Restricted Isometry Property}
In compressed sensing and matrix rank minimization, the restricted
isometry property (RIP) of the matrix $A$ or linear operator $\A$
plays a key role in the relationship between the original
combinatorial problem and its convex relaxation and their optimal
solutions.

The definition of the RIP for matrix rank minimization is:
\begin{definition}\label{def:RIP}
For every integer $r$ with $1\leq r\leq m$, the linear operator $\A:\mathbb{R}^\mtn\rightarrow\mathbb{R}^p$ is
said to satisfy the Restricted Isometry Property with the
restricted isometry constant $\delta_r(\A)$ if $\delta_r(\A)$ is the
minimum constant that satisfies
\bea\label{RIP-definition-eq}(1-\delta_r(\A))\|X\|_F^2\leq\|\A X \|_2^2\leq(1+\delta_r(\A))\|X\|_F^2,\eea
for all $X\in\mathbb{R}^\mtn$ with $\rank(X)\leq r.$ $\delta_r(\A)$
is called the RIP constant. Note that $\delta_s\leq\delta_t,$ if
$s\leq t.$
\end{definition}

The RIP concept and the RIP constant $\delta_r(\A)$ play a central role in the theoretical developments of this paper. We first note that if the operator $\A$ has a nontrivial kernel, i.e., there exists $X\in\br^\mtn$ such that $\A X = 0$ and $X\neq 0$, then $\delta_n(\A) \geq 1$.
Second, if we represent $\A$ in the coordinate form $(\A X)_i=\Tr(A_iX), i=1,\ldots,p,$ then $\delta_r(\A)$ is related to the joint kernel of the matrices $A_i$. For example, if there exists a matrix $X\in\br^\mtn$ with rank $r$ such that $A_iX=0,i=1,\ldots,p,$ then $\delta_r(\A) \geq 1$. Our results in this paper do not apply to such a pathological case.

For matrix rank minimization \eqref{MRM-affine-rank}, Recht \etal
\cite{Recht-Fazel-Parrilo-2007} proved the following results.

\begin{theorem}[Theorem 3.3 in \cite{Recht-Fazel-Parrilo-2007}]\label{thm:Recht-3.3}
Suppose that $\rank(X)\leq r$, $r\geq 1$ and $\delta_{5r}(\A)<0.1$. Then \eqref{MRM-affine-rank} and \eqref{MRM-affine-nuclear} have the same optimal solution.
\end{theorem}

\begin{theorem}[Theorem 4.2 in
\cite{Recht-Fazel-Parrilo-2007}]\label{thm:Recht-4.2} Fix
$\delta\in(0,1)$. If $\A:\br^\mtn \rightarrow \br^p$ is a nearly isometric random map (see
Definition 4.1 in \cite{Recht-Fazel-Parrilo-2007}), then for every
$1\leq r\leq m$, there exist constants $c_0,c_1>0$ depending
only on $\delta$ such that, with probability at least
$1-\exp(-c_1p)$, $\delta_r(\A)\leq \delta$ whenever $p\geq
c_0r(m+n)\log(mn)$.
\end{theorem}

Theorems \ref{thm:Recht-3.3} and \ref{thm:Recht-4.2} indicate that if $\A$ is a nearly
isometric random map, then with very high probability, $\A$ will
satisfy the RIP with a small RIP constant and thus we can solve
\eqref{MRM-affine-rank} by solving its convex relaxation
\eqref{MRM-affine-nuclear}. For example, if $A$ is the matrix
version of the operator $\A$, and its
entries $A_{ij}$ are independent, identically distributed (i.i.d.)
Gaussian, i.e., $A_{ij}\sim\mathcal{N}(0,1/p)$, then $\A$ is a
nearly isometric random map.
For other nearly isometric random maps, see
\cite{Recht-Fazel-Parrilo-2007}.

In Section \ref{sec:numerical-results}, we will show empirically
that when the entries of $A$ are i.i.d. Gaussian, the algorithms
proposed in this paper can solve the matrix rank minimization
problem \eqref{MRM-affine-rank} very well.

It is worth noticing that the linear map $\A$ in the matrix
completion problem \eqref{MRM-MatCompletion-rank} does not satisfy
the RIP. A counterexample is given in \cite{Candes-Plan-2009}. For
more theory on and algorithms for the matrix completion problem, see
\cite{Candes-Recht-2008,Candes-Tao-2009,Candes-Plan-2009,Keshavan-Montanari-Oh-2009-1,Keshavan-Montanari-Oh-2009-2,Cai-Candes-Shen-2008,Ma-Goldfarb-Chen-2008,
Toh-Yun-2009,Liu-Sun-Toh-2009}.

In our proofs of the convergence of FPCA variants, we need
$\A$ to satisfy the RIP. Before we describe some properties of the
RIP that we will use in our proofs, we need the following definitions.

\begin{definition}[{\bf Orthonormal basis of a subspace}]\label{def:basis of subspace}
Given a set of rank-one matrices $\Psi = \{\psi_1,\ldots,\psi_r\}$,
there exists a set of orthonormal matrices $\Gamma =
\{\gamma_1,\ldots,\gamma_s\}$, i.e., $\langle \gamma_i,\gamma_j \rangle =0,$ for
$i\neq j$ and $\|\gamma_i\|_F=1$ for all $i$, such that $\myspan
(\Gamma) =\myspan(\Psi)$. We call $\Gamma$ an {\bf orthonormal
basis} for the subspace $\myspan(\Psi)$. We use $P_\Gamma X$ to
denote the projection of $X$ onto the subspace $\myspan(\Gamma).$
Note that $P_\Gamma X = P_\Psi X$ and $\rank(P_\Gamma X)\leq r,
\forall X\in\mathbb{R}^\mtn.$
\end{definition}

\begin{definition}[{\bf SVD basis of a matrix}]\label{def:basis of rank r matrix}
Assume that the rank-$r$ matrix $X_r$ has the singular value
decomposition $X_r=\sum_{i=1}^r\sigma_iu_iv_i^\top$. $\Gamma :=
\{u_1v_1^\top,u_2v_2^\top,\ldots,u_rv_r^\top\}$ is called an {\bf
SVD basis} for the matrix $X_r.$ Note that elements in $\Gamma$ are
orthonormal rank-one matrices.
\end{definition}

We now list some important properties of linear operators that
satisfy RIP. \footnote{Propositions \ref{prop:RIP:PA*,PA*AP} and
\ref{prop:RIP:error bound of AX} were first proposed by Lee and Bresler without proof in
\cite{Lee-Bresler-2009}. Proofs of Propositions \ref{prop:RIP:PA*,PA*AP} and
\ref{prop:RIP:error bound of AX} were provided later in \cite{Lee-Bresler-ADMIRA-fullpaper-2009}.}

\begin{proposition}\label{prop:RIP:PA*,PA*AP}
Suppose that the linear operator
$\A:\mathbb{R}^\mtn\rightarrow\mathbb{R}^p$ satisfies the RIP with
constant $\delta_r(\A)$. Let $\Psi$ be an arbitrary orthonormal
subset of $\mathbb{R}^\mtn$ such that $\rank(P_\Psi X)\leq r,
\forall X \in \mathbb{R}^\mtn$. Then, for all
$b\in\mathbb{R}^p$ and
$X\in\mathbb{R}^\mtn$, the following properties hold:
\begin{gather}
\|P_\Psi\A^* b\|_F\leq\sqrt{1+\delta_r(\A)}\|b\|_2 \label{prop:RIP:PA*,PA*AP-eq1}\\
(1-\delta_r(\A))\|P_\Psi X\|_F \leq \|P_\Psi\A^*\A P_\Psi X\|_F\leq
(1+\delta_r(\A))\|P_\Psi X\|_F. \label{prop:RIP:PA*,PA*AP-eq2}
\end{gather}
\end{proposition}

\begin{proposition}\label{prop:RIP:PPsi1A*APPsi2}
Suppose that the linear operator
$\A:\mathbb{R}^\mtn\rightarrow\mathbb{R}^p$ satisfies the RIP with
constant $\delta_r(\A)$. Let $\Psi,\Psi'$ be arbitrary orthonormal
subsets of $\mathbb{R}^\mtn$ such that $\rank(P_{\Psi\cup\Psi'}
X)\leq r$, for any $X\in\mathbb{R}^\mtn$. Then the following
inequality holds
\begin{align}\label{prop:RIP:PPsi1A*APPsi2-result}
\|P_\Psi\A^*\A(I-P_\Psi)X\|_F \leq \delta_{r}(\A)\|(I-P_\Psi) X\|_F,
\forall X\in\myspan(\Psi').
\end{align}

\end{proposition}

\begin{proposition}\label{prop:RIP:error bound of AX}
If a linear map $\A:\mathbb{R}^\mtn\rightarrow\mathbb{R}^p$
satisfies \begin{equation}\label{prop:RIP:error bound of
AX-hypothesis-eq}\|\A X\|_2^2\leq
(1+\delta_r(\A))\|X\|_F^2,\quad\forall X\in\mathbb{R}^\mtn,
\rank(X)\leq r,\end{equation} then
\begin{equation}\label{prop:RIP:error bound of AX-conclusion-eq}\|\A X\|_2\leq
\sqrt{1+\delta_r(\A)}\left(\|X\|_F+\frac{1}{\sqrt{r}}\|X\|_*\right),\quad\forall
X\in\mathbb{R}^\mtn.\end{equation}
\end{proposition}

Proofs of Propositions \ref{prop:RIP:PA*,PA*AP},
\ref{prop:RIP:PPsi1A*APPsi2} and \ref{prop:RIP:error bound of AX}
are given in the Appendix.

\section{FPC Revisited} To describe FPC and FPCA and its variants, we need the following definitions.

\begin{definition}\label{def:best rank r approximation}
Assume that the singular value decomposition of the matrix
$X\in\mathbb{R}^\mtn$ is given by
$X=\sum_{i=1}^{m}\sigma_iu_iv_i^\top$ with $\sigma_1\geq\sigma_2\geq\ldots\geq\sigma_m$. Then the best
rank-$r$ approximation $R_r(X)$ to the matrix $X$ is defined as
\begin{align*}R_r(X)=\sum_{i=1}^r\sigma_iu_iv_i^\top.\end{align*}
$R_r: \br^\mtn\rightarrow\br^\mtn$ is also called the hard
thresholding/shrinkage operator with threshold $r$. 
\end{definition}

\begin{definition}\label{def:matrix shrinkage operator}
Assume the SVD of the matrix $X$ is given by
$X=U\Diag(\sigma)V^\top$. For $\nu>0$, the matrix shrinkage operator
$S_\nu(X)$ is defined as
\begin{align*} S_\nu(X) = U \Diag((\sigma-\nu)_+)V^\top, \end{align*}
where $a_+:=\max(a,0)$. $S_\nu:\br^\mtn\rightarrow\br^\mtn$ is also called the soft shrinkage operator with threshold
$\nu$.
\end{definition}

FPC, whose development was motivated by the work on $\ell_1$
regularized problems in \cite{Hale-Yin-Zhang-SIAM-2008}, is based on
applying an operator splitting technique to the optimality
conditions for \eqref{MRM-affine-unconstrained}. Note that $X^*$ is
the optimal solution to \eqref{MRM-affine-unconstrained} if and only
if
\bea\label{OPTcond:MRM-affine-unconstrained}\zero\in\mu\partial\|X^*\|_*+g(X^*),\eea
where $g(X^*)=\A^*(\A X^*-b)$ is the gradient of the least squares
term $\half\|\A X^*-b\|_2^2$, and $\partial\|X^*\|_*$ is the
subgradient of the nuclear norm $\|X^*\|_*$ of $X^*$. According to
\cite{Borwein-Lewis-2000-book}, the subgradient of $\|X\|_*$ is
given by
\begin{align}\label{subgradient-nuclear-norm}\partial\|X\|_*=\{UV^\top+W:U^\top W=0,WV=0,\|W\|_2\leq
1\}, \end{align} where the SVD of $X$ is given by
$X=U\Diag(\sigma)V^\top, U\in\mathbb{R}^{m\times r},
V\in\mathbb{R}^{n\times r},\sigma\in\mathbb{R}_+^r$.

Based on the optimality conditions
\eqref{OPTcond:MRM-affine-unconstrained}, we can develop a fixed-point iterative scheme for solving \eqref{MRM-affine-unconstrained}
by adopting an operator splitting technique. Note that
\eqref{OPTcond:MRM-affine-unconstrained} is equivalent to
\bea\label{optcond:split}\zero\in\tau\mu\partial\|X^*\|_*
+X^*-(X^*-\tau g(X^*))\eea for any $\tau > 0$. If we let \bean Y^* =
X^*-\tau g(X^*),\eea then \eqref{optcond:split} can be rewritten as
\bea\label{optcond:reduced}\zero\in\tau\mu\partial\|X^*\|_*
+X^*-Y^*,\eea i.e., $X^*$ is the optimal solution to
\bea\label{prob:unique-solution-shrinkage-Y*}\min_{X\in\br^{m\times
n}}\tau\mu\|X\|_*+\df{1}{2}\|X-Y^*\|_F^2.\eea

It is known that $S_{\tau\mu}(Y^*)$ gives the optimal solution to
\eqref{prob:unique-solution-shrinkage-Y*}
\cite{Ma-Goldfarb-Chen-2008}. Hence, the following fixed-point
iterative scheme can be given for solving
\eqref{MRM-affine-unconstrained}:
\bea\label{eq:fpc-one-step-scheme}\left\{\ba{l}Y^{k+1}=X^k-\tau
g(X^k)\\X^{k+1}=S_{\tau\mu}(Y^{k+1}).\ea\right.\eea

The following convergence result is proved in
\cite{Ma-Goldfarb-Chen-2008}.
\begin{theorem}[Theorem 4 in \cite{Ma-Goldfarb-Chen-2008}] Assume
$\tau\in(0,2/\lambda_{max}(\A^*\A))$, where $\lambda_{max}(\A^*\A))$ denotes the largest eigenvalue of $\A^*\A$. The sequence $\{X^{k}\}$ generated by the fixed-point
iterations \eqref{eq:fpc-one-step-scheme} converges to some
$X^*\in\X^*,$ where $\X^*$ is the optimal set of problem
\eqref{MRM-affine-unconstrained}.
\end{theorem}

Note that in every iteration of \eqref{eq:fpc-one-step-scheme}, an SVD has to be computed to perform
the matrix shrinkage operation, which is very expensive. Consequently, FPCA uses an approximate SVD to replace the whole SVD, i.e.,
it computes only a rank-$r$ approximation to $Y^{k+1}$. Note that there are
many ways to get a rank-r approximation to $Y^{k+1}$. Here we assume
that the best rank-r approximation $R_r(Y^{k+1})$ is used. In
Section \ref{sec:practical-issues}, we discuss a Monte Carlo
method to approximately compute $R_r(Y^{k+1})$, since computing
$R_r(Y^{k+1})$ exactly is still expensive if $r$ is not very small
and the matrices are large. By adopting a continuation strategy for
the parameter $\mu$ in \eqref{eq:fpc-one-step-scheme}, we arrive at
the following FPCA algorithm (Algorithm \ref{alg:FPCA}) as proposed in \cite{Ma-Goldfarb-Chen-2008}.

\begin{algorithm2e}\caption{Fixed-Point Continuation with Approximate SVD for MRM (FPCA)}
       \label{alg:FPCA}
\linesnumberedhidden \dontprintsemicolon \SetKwInOut{Input}{Input}
\SetKwInOut{Output}{Output} \textsf{Initialization:}  Set $X :=
X^0$. \;
 \For{$\mu=\mu_1,\mu_2,\ldots,\mu_L=\bar\mu$}{
 \While{not converged}{
 $Y := X-\tau\A^*(\A X-b)$. \;
 choose $r$. \;
 $X := S_{\tau\mu}(R_r(Y))$. \;
 }
}
\end{algorithm2e}

We can see that FPCA makes use of three techniques, hard thresholding,
soft shrinkage and continuation. These three techniques have
different properties which, when combined, produce a very robust and
efficient algorithm with great recoverability properties. By using
only one or two of these three techniques, we get different variants
of FPCA. We will study two of these variants, Iterative Hard
Thresholding (IHT) and Iterative Hard Thresholding with soft Matrix
Shrinkage (IHTMS) in Sections \ref{sec:IHT} and \ref{sec:IHTMS},
respectively, and FPCA with given rank $r$ (FPCAr) in Section \ref{sec:fpca}.

In the following three sections, we assume that the rank $r$ of the
optimal solution is given and we compute the best rank-$r$
approximation to $Y$ in each iteration. In Section
\ref{sec:practical-issues}, we give a heuristic for choosing $r$ in
each iteration if $r$ is unknown and use the fast Monte Carlo
algorithm proposed in \cite{Drineas-Kannan-Mahoney-2006} to compute
a rank-$r$ approximation to $Y$.

\section {Iterative Hard Thresholding}\label{sec:IHT} In this section, we study a variant of FPCA that we call Iterative Hard Thresholding (IHT)
because of its similarity to the algorithm in
\cite{Blumensath-Davies-2009} for compressed sensing.

If in FPCA, we assume that the rank $r$ is given, we do not do any
continuation or soft shrinkage, and always choose the stepsize
$\tau$ equal to one, then FPCA becomes Algorithm \ref{alg:IHT}
(IHT). At each iteration of IHT, we first perform a gradient step
$Y^{k+1} := X^k-\A^*(\A X^k-b)$, and then apply hard thresholding
to the singular values of $Y^{k+1}$, i.e., we only keep the largest
$r$ singular values of $Y^{k+1}$, to get $X^{k+1}$.

\begin{algorithm2e}\caption{Iterative Hard Thresholding (IHT)} \label{alg:IHT}

\linesnumberedhidden \dontprintsemicolon \SetKwInOut{Input}{Input}
\SetKwInOut{Output}{Output} \textsf{Initialization:} Given $X^0, r.$
\;
 \For {k = 0,1,\ldots}{
 $Y^{k+1} := X^k-\A^*(\A X^k-b)$. \;
 $X^{k+1}:= R_r(Y^{k+1})$ \;
}
\end{algorithm2e}

As previously mentioned, IHT is closely related to an algorithm
proposed by Blumensath and Davies \cite{Blumensath-Davies-2009} for
compressed sensing. Their algorithm for solving
\eqref{compressed-sensing-L0} performs the following iterative
scheme:

\bea\label{alg:IHT-compressed-sensing}\left\{\ba{l}y^{k+1}=x^k-\tau
A^\top(Ax^k-b) \\x^{k+1}=H_{r}(y^{k+1}),\ea\right.\eea where
$H_r(y)$ is the hard thresholding operator
that sets all but the largest (in magnitude) $r$ elements of $y$ to
zero. Clearly, IHT for matrix rank minimization and compressed
sensing are the same except that the shrinkage operator in the
matrix case is applied to the singular values, while in the
compressed sensing case it is applied to the solution vector.

To prove the convergence/recoverability properties of IHT for matrix
rank minimization, we need the following lemma.

\begin{lemma}\label{lemma:best-rank-r-IHT}
Suppose $X:=R_r(Y)$ is the best rank-$r$ approximation to the matrix
$Y$, and $\Gamma$ is an SVD basis of $X$. Then for any rank-$r$
matrix $X_r$ and SVD basis $\Gamma_r$ of $X_r$, we have
\bea\label{lemma:best-rank-r-IHT-inequality} \|P_BX-P_BY\|_F\leq
\|P_BX_r-P_BY\|_F,\eea where $B$ is any orthonormal set of matrices
satisfying $\myspan(\Gamma\cup\Gamma_r)\subseteq\myspan(B)$.
\end{lemma}
\begin{proof}
Since $X$ is the best rank-$r$ approximation to $Y$ and $rank(X_r)=r$, $\|X-Y\|_F\leq \|X_r-Y\|_F.$
Hence,
\begin{align*}\|P_B(X-Y)\|_F^2+\|(I-P_B)(X-Y)\|_F^2\leq \|P_B(X_r-Y)\|_F^2 + \|(I-P_B)(X_r-Y)\|_F^2.\end{align*}
Since $(I-P_B)X=0$ and $(I-P_B)X_r=0$, this reduces to
\eqref{lemma:best-rank-r-IHT-inequality}.
\end{proof}

For IHT, we have the following convergence results, whose proofs
essentially follow those given by Blumensath and Davies
\cite{Blumensath-Davies-2009} for IHT for compressed sensing.
Our first result considers the case where the desired solution $X_r$
satisfies a perturbed linear system of equations $\A X_r+e=b$.

\begin{theorem}\label{thm:FPCA1-Xs}
Suppose that $b=\A X_r+e$, where $X_r$ is a rank-$r$ matrix, and
$\A$ has the RIP with $\delta_{3r}(\A)\leq\alpha/\sqrt{8}$ where $\alpha\in(0,1)$. Then, at
iteration $k$, IHT will recover an approximation $X^k$ satisfying
\bea\label{Thm-1-FPCA1-eq-1}\left\|X_r-X^k\right\|_F\leq
\alpha^k\left\|X_r-X^0\right\|_F+\frac{\beta}{1-\alpha}\|e\|_2,\eea where $\beta:=2\sqrt{1+\alpha/\sqrt{8}}$. Furthermore, after at most $k^*:=\left\lceil\log_{1/\alpha}\left(\left\|X_r-X^0\right\|_F/\|e\|_2\right)\right\rceil$
iterations, IHT estimates $X^r$ with accuracy
\bea\label{Thm-1-FPCA1-accuracy-k*}\left\|X_r-X^{k^*}\right\|_F\leq
\frac{1-\alpha+\beta}{1-\alpha}\|e\|_2.\eea
\end{theorem}
\begin{proof}
Let $\Gamma_r$ and $\Gamma^k$ denote SVD bases of $X_r$ and $X^k$, respectively, and $B_k$ denote an orthonormal basis of the
subspace $\myspan(\Gamma_r\cup\Gamma^k)$. Let $Z^k:=X_r-X^k$ denote the residual at iteration $k$. Since $P_{B_{k+1}}X_r=X_r$ and
$P_{B_{k+1}}X^{k+1}=X^{k+1}$, it follows first from the triangle
inequality and then from Lemma \ref{lemma:best-rank-r-IHT} that
\begin{equation}\label{proof-Thm-1-FPCA1-triangle-1-leq2}\begin{split}\left\|X_r-X^{k+1}\right\|_F & \leq\left\|P_{B_{k+1}}X_r-P_{B_{k+1}}Y^{k+1}\right\|_F+\left\|P_{B_{k+1}}X^{k+1}-P_{B_{k+1}}Y^{k+1}\right\|_F \\
                                                                   & \leq 2\left\|P_{B_{k+1}}X_r-P_{B_{k+1}}Y^{k+1}\right\|_F.\end{split}\end{equation}
Using the fact that
$b=\A X_r+e$, $Y^{k+1}=X^k-\A^*(\A X^k-\A X_r-e)=X^{k}+\A^*(\A Z^k+e).$
Hence, from \eqref{proof-Thm-1-FPCA1-triangle-1-leq2},
\begin{equation*}
\begin{split}
\left\|X_r-X^{k+1}\right\|_F & \leq 2\left\|P_{B_{k+1}}X_r-P_{B_{k+1}}Y^{k+1}\right\|_F \\ & \leq
2\left\|P_{B_{k+1}}X_r-P_{B_{k+1}}X^k-P_{B_{k+1}}\A^*\A (P_{B_{k+1}}Z^k+(I-P_{B_{k+1}})Z^k) - P_{B_{k+1}}\A^* e\right\|_F \\
                    & \leq 2\left\|P_{B_{k+1}}Z^k-P_{B_{k+1}}\A^*\A(P_{B_{k+1}}Z^k+ (I- P_{B_{k+1}})Z^k)\right\|_F +
                    2\left\|P_{B_{k+1}}\A^* e\right\|_F \\
                    & \leq 2\left\|(I-P_{B_{k+1}}\A^*\A P_{B_{k+1}})P_{B_{k+1}}Z^k\right\|_F + 2\left\|P_{B_{k+1}}\A^*\A (I-P_{B_{k+1}})Z^k\right\|_F +
                    2\left\|P_{B_{k+1}}\A^* e\right\|_F.
\end{split}
\end{equation*}
Since $\rank(P_{B_{k+1}}X)\leq 2r,\forall X\in\mathbb{R}^\mtn$, by
applying \eqref{prop:RIP:PA*,PA*AP-eq1} in Proposition
\ref{prop:RIP:PA*,PA*AP} we get,
\begin{align*}
\left\|P_{B_{k+1}}\A^* e\right\|_F \leq \sqrt{1+\delta_{2r}(\A)}\|e\|_2.
\end{align*}
Since $P_\Psi P_\Psi = P_\Psi$, it follows from
\eqref{prop:RIP:PA*,PA*AP-eq2} in Proposition
\ref{prop:RIP:PA*,PA*AP} that the eigenvalues of the linear operator
$P_\Psi\A^*\A P_\Psi$ 
are in the
interval $[1-\delta_r(\A), 1+\delta_r(\A)]$. Letting $\Psi=B_{k+1}$,
it follows that the eigenvalues of $P_{B_{k+1}}\A^*\A P_{B_{k+1}}$
lie in the interval
$[1-\delta_{2r}(\A), 1+\delta_{2r}(\A)]$. Hence the eigenvalues of
$I-P_{B_{k+1}}\A^*\A P_{B_{k+1}}$
are bounded above by $\delta_{2r}(\A)$ and it follows that
\begin{align*}
\left\|(I-P_{B_{k+1}}\A^*\A P_{B_{k+1}})P_{B_{k+1}}Z^k\right\|_F \leq
\delta_{2r}(\A)\left\|P_{B_{k+1}}Z^k\right\|_F.
\end{align*}
Also, since $P_{B_k}Z^k=Z^k$, $Z^k\in\myspan(B_k)$ and
$\rank(P_{B_k\cup B_{k+1}}X)\leq 3r, \forall X\in\mathbb{R}^\mtn$,
by applying Proposition \ref{prop:RIP:PPsi1A*APPsi2} we get
\begin{align*}
\left\|P_{B_{k+1}}\A^*\A (I-P_{B_{k+1}})Z^k\right\|_F \leq
\delta_{3r}(\A)\left\|(I-P_{B_{k+1}})Z^k\right\|_F.
\end{align*}

Thus, since $\delta_{2r}(\A)\leq\delta_{3r}(\A)$,
\begin{equation*}
\begin{split}
\left\|X_r-X^{k+1}\right\|_F   & \leq 2\delta_{2r}(\A)\left\|P_{B_{k+1}}Z^k\right\|_F + 2
                    \delta_{3r}(\A)\left\|(I-P_{B_{k+1}})Z^k\right\|_F+
                    2\sqrt{1+\delta_{2r}(\A)}\|e\|_2 \\
                    & \leq 2\sqrt{2}\delta_{3r}(\A)\left\|Z^k\right\|_F +
                    2\sqrt{1+\delta_{3r}(\A)}\|e\|_2.
\end{split}
\end{equation*}
By assumption, $\delta_{3r}(\A) \leq \alpha/\sqrt{8}$; hence we have
\bea\label{proof-Thm-1-FPCA1-conclusion-1}\left\|Z^{k+1}\right\|_F \leq \alpha
\left\|Z^k\right\|_F + \beta\|e\|_2.\eea Iterating this inequality, we get
\eqref{Thm-1-FPCA1-eq-1}.

From \eqref{Thm-1-FPCA1-eq-1}, the recovery accuracy $\left\|Z^{k}\right\|_F\leq
\frac{1-\alpha+\beta}{1-\alpha}\|e\|_2$, if $\alpha^k\left\|X_r-X^0\right\|_F\leq
\|e\|_2$. Hence for $k^*:=\left\lceil\log_{1/\alpha}\left(\left\|X_r-X^0\right\|_F/\|e\|_2\right)\right\rceil$,
\eqref{Thm-1-FPCA1-accuracy-k*} holds.
\end{proof}

\begin{remark}
Note that in Theorem \ref{thm:FPCA1-Xs}, convergence is guaranteed for any $\alpha\in(0,1)$. For the choice $\alpha = \half$, $\beta=2\sqrt{1+1/\sqrt{32}}\approx 2.1696$. Thus \eqref{Thm-1-FPCA1-eq-1} becomes
\beaa \left\|X_r-X^k\right\|_F\leq
2^{-k}\left\|X_r-X^0\right\|_F+ 4.3392\|e\|_2,\eeaa and \eqref{Thm-1-FPCA1-accuracy-k*} becomes \beaa \left\|X_r-X^{k^*}\right\|_F\leq
5.3392\|e\|_2. \eeaa
\end{remark}

For an arbitrary matrix $X$, we have the following result.
\begin{theorem}\label{thm:FPCA1-X-arbitrary} Suppose that $b=\A X+e$, where $X$ is an arbitrary matrix, and $\A$ has the RIP
with $\delta_{3r}(\A) \leq \alpha/\sqrt{8}$ where $\alpha\in(0,1)$. Let $X_r$ be the best rank-$r$
approximation to $X$. Then, at iteration $k$, IHT will recover an
approximation $X^k$ satisfying
\bea\label{Thm-2-FPCA1-eq-1}\left\|X-X^k\right\|_F\leq
\alpha^k\left\|X_r-X^0\right\|_F+\gamma\tilde{\epsilon}_r,\eea where $\gamma:=\frac{\beta^2}{2(1-\alpha)}+1$, $\beta:=2\sqrt{1+\alpha/\sqrt{8}}$, and
\bea\label{Thm-2-FPCA1-eq-epsilon}\tilde{\epsilon}_r =
\left\|X-X_r\right\|_F+\frac{1}{\sqrt{r}}\left\|X-X_r\right\|_*+\|e\|_2,\eea is called the
unrecoverable energy (see \cite{Needell-Tropp-2008}). Furthermore,
after at most $k^*:=\left\lceil\log_{1/\alpha}\left(\left\|X_r-X^0\right\|_F/\tilde{\epsilon}_r\right)\right\rceil$
iterations, IHT estimates $X$ with accuracy
\bea\label{Thm-2-FPCA1-accuracy-k*}\left\|X-X^{k^*}\right\|_F\leq
(1+\gamma)\tilde{\epsilon}_r.\eea
\end{theorem}

\begin{proof}
From Theorem \ref{thm:FPCA1-Xs} with $\tilde{e}=\A(X-X_r)+e$ instead of
$e$, we have
\begin{align*}\left\|X_r-X^k\right\|_F=\left\|Z^k\right\|_F\leq \alpha^k\left\|X_r-X^0\right\|_F+\frac{\beta}{1-\alpha}\|\tilde{e}\|_2.\end{align*}
By Proposition \ref{prop:RIP:error bound of AX}, we know that
\begin{align*}\|\tilde{e}\|_2\leq\left\|\A(X-X_r)\right\|_F+\|e\|_2\leq\sqrt{1+\delta_r(\A)}\left(\left\|X-X_r\right\|_F+\frac{1}{\sqrt{r}}\left\|X-X_r\right\|_*\right)+\|e\|_2.\end{align*}
Thus we have from the triangle inequality and \eqref{Thm-2-FPCA1-eq-epsilon}
\begin{align*}\begin{split} \left\|X-X^k\right\|_F & \leq \left\|X_r-X^k\right\|_F + \left\|X-X_r\right\|_F  \\ & \leq \alpha^k\left\|X_r-X^0\right\|_F +
\frac{\beta}{1-\alpha}\|\tilde{e}\|_2 + \left\|X-X_r\right\|_F \\
                                        & \leq \alpha^k\left\|X_r-X^0\right\|_F +
                                        \left(\frac{\beta}{1-\alpha}\sqrt{1+\delta_r(\A)}+1\right)\tilde{\epsilon}_r \\
                                        & \leq \alpha^k\left\|X_r-X^0\right\|_F+\gamma \tilde{\epsilon}_r.\end{split}\end{align*}
This proves \eqref{Thm-2-FPCA1-eq-1}.

Furthermore, $\left\|X-X^k\right\|_F\leq (1+\gamma)\tilde{\epsilon}_r$ if $\alpha^k\left\|X_r-X^0\right\|_F\leq\tilde{\epsilon}_r$.
Therefore, for $k^*:=\left\lceil\log_{1/\alpha}\left(\left\|X_r-X^0\right\|_F/\tilde{\epsilon}_r\right)\right\rceil$,
\eqref{Thm-2-FPCA1-accuracy-k*} holds.
\end{proof}

\begin{remark}
For the choice $\alpha = \half$, $\beta=2\sqrt{1+1/\sqrt{32}}\approx 2.1696$ and $\gamma=\frac{\beta^2}{2(1-\alpha)}+1\approx 5.7072$. Thus \eqref{Thm-2-FPCA1-eq-1} holds as \beaa \left\|X-X^k\right\|_F\leq
2^{-k}\left\|X_r-X^0\right\|_F+ 5.7072\tilde{\epsilon}_r,\eeaa and \eqref{Thm-2-FPCA1-accuracy-k*} holds as \beaa \left\|X-X^{k^*}\right\|_F\leq
6.7072\tilde{\epsilon}_r.\eeaa
\end{remark}

Similar bounds on the RIP constant for an approximate recovery were obtained by Lee and Bresler \cite{Lee-Bresler-2009,Lee-Bresler-ADMIRA-fullpaper-2009} for affinely constrained matrix rank minimization and by Lee and Bresler for ellipsoidally constrained matrix rank minimization \cite{Lee-Bresler-Ellipsoidal-2009}.
The results in Theorems \ref{thm:FPCA1-Xs} and \ref{thm:FPCA1-X-arbitrary} improve the previous results for affinely constrained matrix rank minimization in \cite{Lee-Bresler-2009,Lee-Bresler-ADMIRA-fullpaper-2009}. Specifically, Theorems \ref{thm:FPCA1-Xs} and \ref{thm:FPCA1-X-arbitrary} require the RIP constant $\delta_{3r}(\A) < 1/\sqrt{8} \approx 0.3536$, while the result in \cite{Lee-Bresler-2009,Lee-Bresler-ADMIRA-fullpaper-2009} requires $\delta_{4r}(\A) \leq 0.04$ and the result in \cite{Lee-Bresler-Ellipsoidal-2009} requires $\delta_{3r}(\A) < 1/(1+4/\sqrt{3}) \approx 0.3022$ for recovery in the noisy case.
The IHT algorithm for matrix rank minimization has also been independently studied by Meka, Jain and Dhillon in \cite{Meka-Jain-Dhillon-SVP-2009}, who have obtained very different results than those in Theorems \ref{thm:FPCA1-Xs} and \ref{thm:FPCA1-X-arbitrary}.

\section {Iterative Hard Thresholding with Matrix Shrinkage}
\label{sec:IHTMS}
We study another variant of FPCA in this section. If in each
iteration of IHT, we perform matrix shrinkage to $R_r(Y)$ with fixed
thresholding $\mu>0$, we get the following algorithm (Algorithm \ref{alg:IHTMS}), which we call
Iterative Hard Thresholding with Matrix Shrinkage (IHTMS).
Note that $S_{\mu}(R_r(Y))=R_r(S_\mu(Y)),\forall r,\mu$ and $Y$.

\begin{algorithm2e}\caption{Iterative Hard Thresholding with Matrix Shrinkage (IHTMS)} \label{alg:IHTMS}

\linesnumberedhidden \dontprintsemicolon \SetKwInOut{Input}{Input}
\SetKwInOut{Output}{Output} \textsf{Initialization:} Given $X^0,
\mu$ and $r$. \;
 \For {k = 0,1,\ldots}{
 $Y^{k+1} := X^k-\A^*(\A X^k-b)$. \;
 $X^{k+1}:=R_r(S_\mu(Y^{k+1}))$. \;
}
\end{algorithm2e}

For IHTMS, we have the following convergence results.

\begin{theorem}\label{thm:FPCA2-Xs} Suppose that
$b=\A X_r+e$, where $X_r$ is a rank-$r$ matrix, and $\A$ has the
RIP with $\delta_{3r}(\A) \leq \alpha/\sqrt{8}$ where $\alpha\in(0,1)$. Then, at iteration $k$,
IHTMS will recover an approximation $X^k$ satisfying
\bea\label{Thm-1-FPCA2-eq-1}\left\|X_r-X^k\right\|_F\leq
\alpha^k\left\|X_r-X^0\right\|_F+  \frac{1}{1-\alpha}(\beta\|e\|_2 + 2\mu\sqrt{m}),\eea where $\beta:=2\sqrt{1+\alpha/\sqrt{8}}$.
Furthermore, after at most $k^*:=\left\lceil\log_{1/\alpha}\left(\left\|X_r-X^0\right\|_F/(\|e\|_2+2\mu\sqrt{m})\right)\right\rceil$
iterations, IHTMS estimates $X_r$ with accuracy
\bea\label{Thm-1-FPCA2-accuracy-k*}\left\|X_r-X^{k^*}\right\|_F\leq \frac{1-\alpha+\beta}{1-\alpha} \|e\|_2 +
\frac{2-\alpha}{1-\alpha}2\mu\sqrt{m}.\eea
\end{theorem}
\begin{proof}
Using the same notation as in the proof of Theorem \ref{thm:FPCA1-Xs}, we know that $P_{B_{k+1}}X_r=X_r$ and $P_{B_{k+1}}X^{k+1}=X^{k+1}$. Using the
triangle inequality we get,
\begin{equation}\label{proof-Thm-1-IHTMS-eq-1}
\begin{split} \left\|X_r-X^{k+1}\right\|_F \leq &
\left\|P_{B_{k+1}}X_r-P_{B_{k+1}}Y^{k+1}\right\|_F \\
& + \left\|P_{B_{k+1}} X^{k+1}-P_{B_{k+1}}S_\mu(Y^{k+1})\right\|_F\\ & +
\left\|P_{B_{k+1}}S_\mu(Y^{k+1})-P_{B_{k+1}}Y^{k+1}\right\|_F.
\end{split}
\end{equation}
Since $X^{k+1}$ is the best rank-$r$ approximation to
$S_\mu(Y^{k+1})$, by applying Lemma \ref{lemma:best-rank-r-IHT} we get
\begin{equation}\label{proof-Thm-1-IHTMS-eq-2}
\begin{split}
\left\|P_{B_{k+1}}X^{k+1}-P_{B_{k+1}}S_\mu(Y^{k+1})\right\|_F
\leq & \left\|P_{B_{k+1}}X_r - P_{B_{k+1}}S_\mu(Y^{k+1})\right\|_F \\  \leq &
\left\|P_{B_{k+1}}X_r-P_{B_{k+1}}Y^{k+1}\right\|_F  \\  & +
\left\|P_{B_{k+1}}S_\mu(Y^{k+1})-P_{B_{k+1}} Y^{k+1}\right\|_F.
\end{split}
\end{equation}
Therefore, by combining \eqref{proof-Thm-1-IHTMS-eq-1},
\eqref{proof-Thm-1-IHTMS-eq-2} and noticing that
\begin{equation*}
\left\|P_{B_{k+1}}S_\mu(Y^{k+1})-P_{B_{k+1}}Y^{k+1}\right\|_F \leq
\left\|S_\mu (Y^{k+1})-Y^{k+1}\right\|_F \leq \mu\sqrt{m},
\end{equation*}
we have
\begin{equation*}
\left\|X_r-X^{k+1}\right\|_F \leq 2\left\|P_{B_{k+1}} X_r-P_{B_{k+1}}Y^{k+1}\right\|_F +
2\mu\sqrt{m}.
\end{equation*}
Using an argument identical the one below \eqref{proof-Thm-1-FPCA1-triangle-1-leq2} in the proof of Theorem \ref{thm:FPCA1-Xs}, we get
\begin{equation*}
\left\|X_r-X^{k+1}\right\|_F \leq 2\sqrt{2}\delta_{3r}(\A)\left\|Z^k\right\|_F +
                    2\sqrt{1+\delta_{3r}(\A)}\|e\|_2+
2\mu\sqrt{m}.
\end{equation*}
Now since $\delta_{3r}(\A) \leq \alpha/\sqrt{8}$, we have
\begin{equation*}
\left\|Z^{k+1}\right\|_F \leq \alpha\left\|Z^k\right\|_F + \beta\|e\|_2 +
2\mu\sqrt{m},
\end{equation*}
which implies that \eqref{Thm-1-FPCA2-eq-1} holds. Hence \eqref{Thm-1-FPCA2-accuracy-k*} holds if $k^*:=\left\lceil\log_{1/\alpha}\left(\left\|X_r-X^0\right\|_F/(\|e\|_2+2\mu\sqrt{m})\right)\right\rceil$.
\end{proof}

For an arbitrary matrix $X$, we have the following results.
\begin{theorem}\label{thm:FPCA2-X-arbitrary} Suppose that $b=\A X+e$, where $X$ is an arbitrary matrix, and $\A$ has the RIP
with $\delta_{3r}(\A)\leq \alpha/\sqrt{8}$ where $\alpha\in(0,1)$. Let $X_r$ be the best rank-$r$
approximation to $X$. Then, at iteration $k$, IHTMS will recover an
approximation $X^k$ satisfying
\bea\label{Thm-2-FPCA2-eq-1}\left\|X-X^k\right\|_F\leq
\alpha^k\left\|X_r-X^0\right\|_F+\gamma\tilde{\epsilon}_r+\frac{2\mu\sqrt{m}}{1-\alpha},\eea where $\gamma:=\frac{\beta^2}{2(1-\alpha)}+1$, $\beta:=2\sqrt{1+\alpha/\sqrt{8}}$, and
$\tilde{\epsilon}_r$ is defined by \eqref{Thm-2-FPCA1-eq-epsilon}. Furthermore,
after at most $k^*:=\left\lceil\log_{1/\alpha}\left(\left\|X_r-X^0\right\|_F/(\tilde{\epsilon}_r+2\mu\sqrt{m})\right)\right\rceil$
iterations, IHTMS estimates $X$ with accuracy
\bea\label{Thm-2-FPCA2-accuracy-k*}\left\|X-X^{k^*}\right\|_F\leq
(1+\gamma)\tilde{\epsilon}_r+\frac{2-\alpha}{1-\alpha}2\mu\sqrt{m}.\eea
\end{theorem}

\begin{proof}
The proof of \eqref{Thm-2-FPCA2-eq-1} is identical to the proof of \eqref{Thm-2-FPCA1-eq-1} in Theorem \ref{thm:FPCA1-X-arbitrary}, except that \eqref{Thm-1-FPCA2-eq-1} is used instead of \eqref{Thm-1-FPCA1-eq-1}. It also immediately follows from \eqref{Thm-2-FPCA2-eq-1} that \eqref{Thm-2-FPCA2-accuracy-k*} holds for $k^*:=\left\lceil\log_{1/\alpha}\left(\left\|X_r-X^0\right\|_F/(\tilde{\epsilon}_r+2\mu\sqrt{m})\right)\right\rceil$.\end{proof}

\section {FPCA with Given Rank $r$}\label{sec:fpca} In this section, we study the
FPCA when rank $r$ is known and a unit stepsize $\tau=1$ is always chosen. This is equivalent to applying a continuation strategy to
$\mu$ in IHTMS. We call this algorithm FPCAr (see Algorithm \ref{alg:FPCAr} below). The parameter $\eta_\mu$ determines the
rate of reduction of the consecutive $\mu_j$ in continuation, i.e.,
\begin{align}\label{dec:mu}\mu_{j+1}=\max\{\mu_j\eta_\mu,\bar\mu\}, j=1,\ldots,L-1\end{align}

For FPCAr, we have the following convergence results.

\begin{algorithm2e}\caption{FPCA with given rank $r$ (FPCAr)}
       \label{alg:FPCAr}
\linesnumberedhidden \dontprintsemicolon \SetKwInOut{Input}{Input}
\SetKwInOut{Output}{Output} 
\Input{$X_{(1)}^0, r, \mu_1>\mu_2\ldots>\mu_L=\bar\mu.$} \;
 \For {j = 1,\ldots,L}{
 Set $\mu = \mu_j.$ \;
 \For {k = 0,1,\ldots, until convergence}{
 $Y_{(j)}^{k+1} := X_{(j)}^k-\A^*\left(\A X_{(j)}^k-b\right)$. \;
 $X_{(j)}^{k+1}:=S_{\mu}\left(R_r\left(Y_{(j)}^{k+1}\right)\right)$. \;
 }
 Set $X_{(j+1)}^0 = X_{(j)}^{k+1}.$ \;
} \Output {$X^*:=X_{(L+1)}^0.$}
\end{algorithm2e}

\begin{theorem}\label{thm:FPCAr-Xs}
Suppose that $b=\A X_r+e$, where $X_r$ is a rank-$r$ matrix, and
$\A$ has the RIP with $\delta_{3r}(\A) \leq \alpha/\sqrt{8}$ where $\alpha\in(0,1)$. Also, suppose
in FPCAr, after $K_j$ iterations with fixed $\mu=\mu_j$, we obtain a
solution $X_{(j)}^{(K_j)}$ that is then set to the initial point
$X_{(j+1)}^0$ for the next continuation subproblem $\mu=\mu_{j+1}$.
Then FPCAr will recover an approximation $X_{(L)}^{(K_L)}$ that
satisfies
\bea\label{Thm-1-FPCA3-KL-accuracy}\begin{split}\left\|X_r-X_{(L)}^{(K_L)}\right\|_F
\leq & \left(\alpha^{\sum_{j=1}^L K_j}\right)\left\|X_r-X^0\right\|_F+\left(\sum_{j=2}^L
\alpha^{\sum_{l=j}^L K_l} + 1\right)\frac{\beta}{1-\alpha}\|e\|_2 \\ & +
\left(\sum_{j=2}^L\left(\alpha^{\sum_{l=j}^L
K_l}\right)\mu_{j-1}+\mu_L\right)\frac{2\sqrt{m}}{1-\alpha},\end{split}\eea where $\beta:=2\sqrt{1+\alpha/\sqrt{8}}$.
\end{theorem}
\begin{proof}
For $X_{(1)}^{(K_1)}$, which is obtained by setting $\mu=\mu_1$ in
the first $K_1$ iterations, we get from
Theorem \ref{thm:FPCA2-Xs}, that if $\delta_{3r}(\A) \leq \alpha/\sqrt{8}$,
\bea\label{proof-Thm-1-FPCA3-conclusion-2-1st-iter}\left\|X_r-X_{(1)}^{(K_1)}\right\|_F
\leq \alpha^{K_1}\left\|X_r-X^0\right\|_F + \frac{\beta}{1-\alpha}\|e\|_2+\frac{2\mu_1\sqrt{m}}{1-\alpha}.\eea
Then from iteration $K_1+1$ to $K_1+K_2$, we fix $\mu=\mu_2$. Again
by Theorem \ref{thm:FPCA2-Xs}, we get
\bea\label{proof-Thm-1-FPCA3-conclusion-2-2nd-iter}\left\|X_r-X_{(2)}^{(K_2)}\right\|_F
\leq \alpha^{K_2}\left\|X_r-X_{(1)}^{(K_1)}\right\|_F +
\frac{\beta}{1-\alpha}\|e\|_2+\frac{2\mu_2\sqrt{m}}{1-\alpha}.\eea
By substituting \eqref{proof-Thm-1-FPCA3-conclusion-2-1st-iter} into
\eqref{proof-Thm-1-FPCA3-conclusion-2-2nd-iter}, we get \begin{eqnarray*}
\label{proof-Thm-1-FPCA3-conclusion-3-2nd-iter}\begin{split}
\left\|X_r-X_{(2)}^{(K_2)}\right\|_F \leq & \alpha^{(K_1+K_2)}\left\|X_r-X^0\right\|_F +
\left(\alpha^{K_2}+1\right)\frac{\beta}{1-\alpha}\|e\|_2
\\ & +\left(\alpha^{K_2}\mu_1+\mu_2\right)\frac{2\sqrt{m}}{1-\alpha}.\end{split}\end{eqnarray*}
Repeating this procedure we get \eqref{Thm-1-FPCA3-KL-accuracy}. 
\end{proof}

Theorem \ref{thm:FPCAr-Xs} shows that as long as $\mu_L$ is small and $K_L$ is large, the recovery
error will be very small.
For an arbitrary matrix $X$, we have the following convergence result.
\begin{theorem}\label{thm:FPCAr-X-arbitrary}
Suppose that $b=\A X+e$, where $X$ is an arbitrary matrix. Let
$X_r$ be the best rank-$r$ approximation to $X$. With the same
notation and under the same conditions as in Theorem
\ref{thm:FPCAr-Xs}, FPCAr will recover an approximation
$X_{(L)}^{(K_L)}$ that satisfies
\beaa\label{Thm-2-FPCA3-KL-accuracy}\begin{split}\left\|X-X_{(L)}^{(K_L)}\right\|_F
\leq & \left(\alpha^{\sum_{j=1}^L K_j}\right)\left\|X_r-X^0\right\|_F+\left(\left(\sum_{j=2}^L
\alpha^{\sum_{l=j}^L K_l} + 1\right)\gamma+1\right)\tilde{\epsilon}_r \\ & +
\left(\sum_{j=2}^L\left(\alpha^{\sum_{l=j}^L
K_l}\right)\mu_{j-1}+\mu_L\right)\frac{2\sqrt{m}}{1-\alpha},\end{split}\eeaa where $\gamma:=\frac{\beta^2}{2(1-\alpha)}+1$, $\beta:=2\sqrt{1+\alpha/\sqrt{8}}$, and $\tilde{\epsilon}_r$ is defined by \eqref{Thm-2-FPCA1-eq-epsilon}.
\end{theorem}
\begin{proof}
We skip the proof here since it is similar to the proof of Theorem
\ref{thm:FPCA1-X-arbitrary}.
\end{proof}

\section{Practical Issues}\label{sec:practical-issues}
In practice, the rank $r$ of the optimal solution is usually unknown. Thus, in every iteration, we need to determine $r$
appropriately. We propose some heuristics for doing this here. We
start with $r:=r_{\max}$. So
$X^1$ is a rank-$r_{\max}$ matrix. For the $k$-th iteration ( $k\geq
2$ ), $r$ is chosen as the number of singular values of $X^{k-1}$
that are greater than $\epsilon_s\sigma_1^{k-1}$, where
$\sigma_1^{k-1}$ is the largest singular value of $X^{k-1}$ and
$\epsilon_s\in(0,1)$ is a given tolerance. Sometimes the given
tolerance truncates too many of the singular values, so we need to
increase $r$ occasionally. One way to do this is to increase $r$ by
1 whenever the non-expansive property (see \cite{Ma-Goldfarb-Chen-2008}) of the shrinkage operator
$S_\mu$ is violated some fixed number of times, say 10.
In the numerical experiments described in Section
\ref{sec:numerical-results}, we used another strategy; i.e., we
increased $r$ by 1 whenever the Frobenius norm of the gradient $g$
increased by more than 10 times. We tested this heuristic for
determining $r$ extensively. It enables our algorithms to achieve
very good recoverability and appears to be very robust. For many
examples, our algorithms can recover matrices whose rank is almost
$r_{\max}$ with a limited number of measurements.

Another issue in practice is concerned with the SVD computation.
Note that in IHT, IHTMS and FPCA, we need to compute the best
rank-$r$ approximation to $Y^{k+1}$ at every iteration. This can be
very expensive even if we use a state-of-the-art code like PROPACK
\cite{Larsen-Propack}, especially when the rank of the matrix is
relatively large. Therefore, we used instead the Monte Carlo
algorithm LinearTimeSVD proposed in
\cite{Drineas-Kannan-Mahoney-2006} to approximate the best rank-$r$
approximation. 
For a given matrix $A\in\mathbb{R}^{m\times n}$, and parameters
$c_s,k_s\in\mathbb{Z}^+$ with $1\leq k_s\leq c_s\leq n$ and
$\{p_i\}_{i=1}^n,$ $p_i\geq 0, \sum_{i=1}^np_i=1$, this algorithm
returns approximations $\sigma_t(C),t=1,\ldots,k_s$ to the largest $k_s$ singular values and approximations $H_{k_s}^{(t)},t=1,\ldots,k$ to the
corresponding left singular vectors of the matrix $A$ in $O(m+n)$
time. Thus, the SVD of $A$ is approximated by $$A\approx
A_{k_s}:= H_{k_s}\Diag(\sigma(C))(A^\top
H_{k_s}\Diag(1/\sigma(C))^\top.$$
Drineas
\etal\cite{Drineas-Kannan-Mahoney-2006} prove that with high
probability, the following estimate holds for both $\xi=2$ and
$\xi=F$ when $\{p_i\}_{i=1}^n$ are {\it nearly optimal probabilities} (see \cite{Drineas-Kannan-Mahoney-2006}):
\bea\label{eq:error-fast-svd}\|A-A_{k_s}\|_\xi^2\leq\min_{D:\rank(D)\leq
k_s}\|A-D\|_\xi^2+poly(k_s,1/c_s)\|A\|_F^2,\eea where
$poly(k_s,1/c_s)$ is a polynomial in $k_s$ and $1/c_s$. Thus,
$A_{k_s}$ is an approximation to the best rank-$k_s$ approximation to
$A$. The LinearTimeSVD Algorithm, which we found to be much faster
than PROPACK, is outlined below in Algorithm \ref{alg:Linear-Time-SVD}.

\begin{algorithm2e}\caption{Linear Time Approximate SVD
Algorithm \cite{Drineas-Kannan-Mahoney-2006}}\label{alg:Linear-Time-SVD}
\linesnumberedhidden \dontprintsemicolon \SetKwInOut{Input}{Input}
\SetKwInOut{Output}{Output} 
\Input{$A\in\mathbb{R}^{m\times n}$, $c_s,k_s\in\mathbb{Z}^+$ \st $1\leq k_s\leq c_s\leq
n$, $\{p_i\}_{i=1}^n$ \st $p_i\geq 0, \sum_{i=1}^np_i=1$.} \;
\Output{$H_k\in\mathbb{R}^{m\times k_s}$ and $\sigma_t(C), t=1,\ldots,k_s.$}
 \For {t = 1,\ldots,$c_s$}{
 Pick $i_t\in 1,\ldots,n$ with $Pr[i_t=\alpha] =
p_\alpha,\alpha=1,\ldots,n.$ \;
Set $C^{(t)}=A^{(i_t)}/\sqrt{c_sp_{i_t}}.$ \;
}
Compute $C^\top C$ and its SVD; say $
C^\top C=\sum_{t=1}^{c_s}\sigma_t^2(C)y^{t}{y^t}^\top.$ \;
Compute $h^t=Cy^t/\sigma_t(C)$ for $t=1,\ldots,k_s.$ \;
Return $H_{k_s}$, where $H_{k_s}^{(t)}=h^t,$ and $\sigma_t(C),
t=1,\ldots,k_s.$ \;
\end{algorithm2e}

Note that in Algorithm \ref{alg:Linear-Time-SVD}, we compute an exact SVD of a smaller
matrix $C^\top C\in\mathbb{R}^{c_s\times c_s}$. Thus, $c_s$
determines the speed of this algorithm. If we choose a large $c_s$,
we need more time to compute the SVD of $C^\top C$. However, the
larger $c_s$ is, the more likely are the
$\sigma_t(C),t=1,\ldots,k_s$ to be close to the largest $k_s$
singular values of the matrix $A$ since the second term in the right
hand side of \eqref{eq:error-fast-svd} is smaller. In our numerical
experiments, we found that we could choose a relatively small $c_s$
so that the computational time was reduced without significantly
degrading the accuracy. There are many ways to choose the probabilities $p_i$. In our
numerical experiments in Section \ref{sec:numerical-results}, we used the simplest one, i.e.,
we set all $p_i$ equal to $1/n$. For other choices of $p_i$, see
\cite{Drineas-Kannan-Mahoney-2006} and the references therein.

Although
PROPACK is more accurate than this Monte Carlo method (Algorithm \ref{alg:Linear-Time-SVD}), we observed from our numerical experiments that our algorithms are very robust and are not very sensitive to the accuracy
of the approximate SVDs.

In the $j$-th inner iteration in FPCA we solve problem
\eqref{MRM-affine-unconstrained} for a fixed $\mu=\mu_j$; and stop
when
\bea\label{stop:fpca}\frac{\|X^{k+1}-X^k\|_F}{\max\{1,\|X^k\|_F\}}<xtol,\eea
where $xtol$ is a small positive number. We then decrease $\mu$ and
go to the next inner iteration.


\section {Numerical Experiments} \label{sec:numerical-results}
In this section, we present numerical results for the algorithms
discussed above and provide comparisons with the SDP solver SDPT3
\cite{Toh-Todd-Tutuncu-SDPT3}. We use IHTr, IHTMSr, FPCAr to denote algorithms in which the rank $r$ is specified, and IHT, IHTMS, FPCA to denote those in which $r$ is determined by the heuristics described in
Section \ref{sec:practical-issues}. We tested these six algorithms on both randomly created
and realistic matrix rank minimization problems \eqref{MRM-affine-rank}. IHTr, IHT, IHTMSr and IHTMS were terminated when \eqref{stop:fpca} holds. FPCAr and FPCA were terminated when both \eqref{stop:fpca} holds and $\mu_{k}=\bar\mu$. All numerical experiments were run in MATLAB 7.3.0 on a Dell Precision 670
workstation with an Intel xeon(TM) 3.4GHZ CPU and 6GB of RAM. All CPU times reported in this section are in seconds.

\subsection{Randomly Created Test Problems}\label{sec:numerical-results-random}
We tested some randomly created problems to illustrate the recoverability/convergence properties of our algorithms. The random test problems \eqref{MRM-affine-rank} were created in the
following manner. We first generated random matrices
$M_L\in\br^{m\times r}$ and $M_R\in\br^{n\times r}$ with i.i.d.
Gaussian entries $\sim\mathcal{N}(0,1)$ and then set
$M=M_LM_R^\top$. We then created a matrix $A\in\br^{p\times mn}$
with i.i.d. Gaussian entries $A_{ij}\sim\mathcal{N}(0,1/p).$
Finally, the observation $b$ was set equal to $b:=A\vvec(M)$. We use
$SR=p/(mn)$, i.e., the number of measurements divided by the number
of entries of the matrix, to denote the sampling ratio. We also list
$FR=r(m+n-r)/p$, i.e. the dimension of the set of rank $r$ matrices
divided by the number of measurements, in the tables. Note that if
$FR>1$, then there is always an infinite number of matrices with
rank $r$ satisfying the $p$ linear constraints, so we cannot hope to
recover the matrix in this situation. We also report the relative
error \bean rel.err. := \frac{\|X_{opt}-M\|_F}{\|M\|_F}\eea to
indicate the closeness of $X_{opt}$ to $M$, where $X_{opt}$ is the
optimal solution to \eqref{MRM-affine-rank} produced by our
algorithms. We declared $M$ to be recovered if the relative error
was less than $10^{-3}.$ We solved 10 randomly created matrix rank
minimization problems for each set of $(m,n,p,r)$. We used $NS$ to
denote the number of matrices that were recovered successfully. The
average time and average relative error of the successfully solved
problems are also reported.

The parameters used in the algorithms are summarized in Table
\ref{table:parameters}.
\begin{table}[ht]
\begin{center}\caption{Parameters used in the algorithms}\label{table:parameters}
\begin{tabular}{c|c|c}\hline
parameter & value & description
\\\hline\hline

$\bar\mu$ & $10^{-8}$ & parameter in Algorithms \ref{alg:FPCA} and
\ref{alg:FPCAr}
\\\hline

$\eta_\mu$ & 0.25 & parameter in \eqref{dec:mu}  \\\hline


$\epsilon_s$ & 0.01 & parameter in LinearTimeSVD \\\hline

$c_s$ & $2r_{\max}-2$ & parameter in LinearTimeSVD \\\hline

$p_i$ & $1/n, \forall i$ & parameter in LinearTimeSVD \\\hline

$xtol$ & $10^{-6}$ & parameter in \eqref{stop:fpca} \\\hline

\end{tabular}
\end{center}
\end{table}

We first compare the solvers discussed above that specify the rank $r$
with the SDP solver SDPT3 \cite{Toh-Todd-Tutuncu-SDPT3}. The
results for a set of small problems with $m=n=60$, 20 percent
sampling (i.e., SR = 0.2 and p = 720) and different ranks are presented
in Table \ref{table:Num-res-all-SDPT3-Given-rank}. Note that for
this set of parameters $(m,n,p)$, the largest rank that satisfies
$FR<1$ is $r_{\max}=6$.
\begin{table}[ht]
\begin{center}\caption{Comparison between IHTr, IHTMSr and FPCAr with SDPT3}\label{table:Num-res-all-SDPT3-Given-rank}
\begin{tabular}{c|c|c|c|c|c|c|c|c|c|c|c|c|c}\hline
\multicolumn{2}{c|}{Prob} & \multicolumn{3}{c|}{SDPT3} &
\multicolumn{3}{c|}{IHTr} & \multicolumn{3}{c|}{IHTMSr} &
\multicolumn{3}{c}{FPCAr} \\\hline

r &  FR  & NS & time & rel.err. & NS & time & rel.err. & NS & time &
rel.err. & NS & time & rel.err.
\\\hline\hline

1    & 0.17 & 10 & 122.93 & 2.31e-10 &10     & 2.60      & 1.67e-05
&10 & 2.59 & 1.67e-05     &10     & 4.63      & 9.00e-06 \\\hline

2    & 0.33 & 10 & 124.26 & 3.46e-09  &10     & 4.97      & 1.99e-05
&10 & 4.98      & 2.11e-05     &10     & 6.06      & 1.51e-05
\\\hline

3    & 0.49 & 3 & 149.74 & 2.84e-07  &10     & 10.04     & 2.38e-05
&10 & 9.95      & 2.27e-05     &10     & 10.64     & 2.35e-05
\\\hline

4    & 0.64 & 0 & --- & --- &10     & 22.99     & 2.88e-05 &10 &
22.72     & 3.05e-05     &10     & 23.29     & 2.93e-05
\\\hline

5    & 0.80 & 0 & --- & --- &10     & 75.86     & 3.89e-05 &10 &
84.13     & 3.95e-05     &10     & 79.46     & 3.94e-05
\\\hline

\end{tabular}
\end{center}
\end{table}

From Table \ref{table:Num-res-all-SDPT3-Given-rank} we can see that
the performance of our methods is very robust and quite similar in terms of their recoverability properties. They are also much faster
and their abilities to recover the matrices are much better than
SDPT3. For ranks less than or equal to 5, which is almost the
largest rank guaranteeing $FR<1$, IHTr, IHTMSr and FPCAr can recover
all randomly generated matrices with a relative error of the order
of $1e-5$. However, SDPT3 can only recover all matrices with a rank
equal to 1 or 2. When the rank $r$ increases to 3, SDPT3 can only
recover 3 of the 10 matrices. When the rank $r$ increases to 4 or 5,
none of the 10 matrices can be recovered by SDPT3.

\begin{figure}
\begin{center}
\includegraphics[scale=0.65]{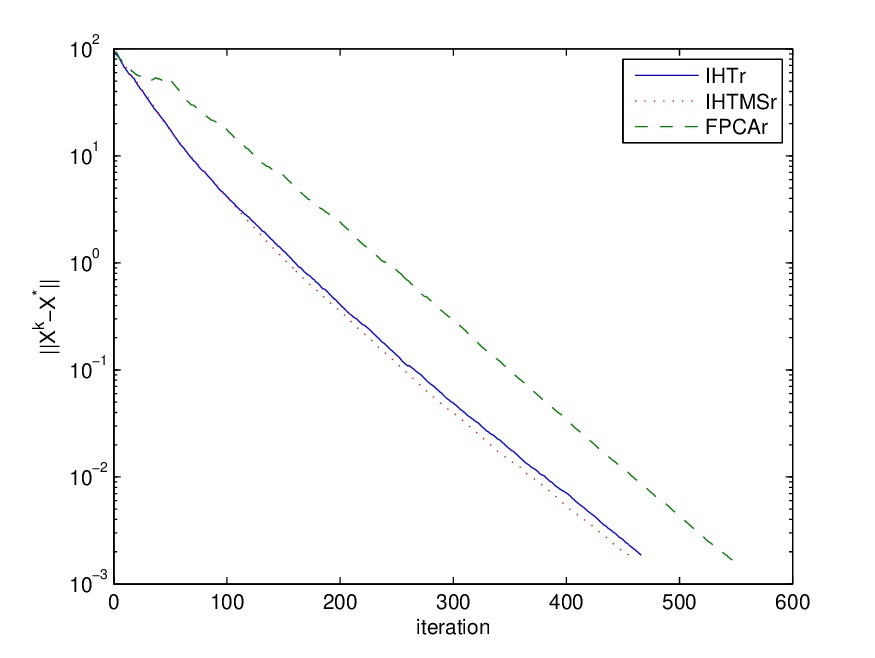}
\caption{Approximation error versus the iteration number for a
problem where the rank equaled 2} \label{fig:error}
\end{center}
\end{figure}

To verify the theoretical results in Sections 4, 5 and 6, we plotted
the log of the approximation error $\|X^k - X^*\|_F$ achieved by each of the algorithms IHTr, IHTMSr and FPCAr versus the
iteration number $k$ in Figure
\ref{fig:error} for one of 10 randomly created problems involving a matrix of rank 2.
From this figure, we can see that $\log\|X^k-X^*\|_F$ is
approximately a linear function of the iteration number $k$. This
implies that our theoretical results in Sections 4, 5 and 6
approximately hold in practice.

For the same set of test problems, Tables \ref{table:Num-res-all-comparison-r-IHT},
\ref{table:Num-res-all-comparison-r-IHTMS}, and
\ref{table:Num-res-all-comparison-r-FPCA} present comparisons of IHTr
versus IHT, IHTMSr versus IHTMS and FPCAr versus FPCA.

\begin{table}[ht]
\begin{center}\caption{Comparison between IHTr and IHT}\label{table:Num-res-all-comparison-r-IHT}
\begin{tabular}{c|c|c|c|c|c|c|c}\hline
\multicolumn{2}{c|}{Prob} & \multicolumn{3}{c|}{IHTr} &
\multicolumn{3}{c}{IHT} \\\hline

r &  FR  & NS & time & rel.err. & NS & time & rel.err.
\\\hline\hline
1    & 0.17      &10     & 2.60      & 1.67e-05         &10     &
4.24      & 1.74e-05 \\\hline

2    & 0.33      &10     & 4.97      & 1.99e-05         &10     &
7.00      & 1.92e-05 \\\hline

3    & 0.49      &10     & 10.04     & 2.38e-05         &10     &
13.27     & 2.32e-05 \\\hline

4    & 0.64      &10     & 22.99     & 2.88e-05         &10     &
28.06     & 2.93e-05 \\\hline

5    & 0.80      &10     & 75.86     & 3.89e-05         &10     &
96.32     & 4.00e-05 \\\hline
\end{tabular}
\end{center}
\end{table}

\begin{table}[ht]
\begin{center}\caption{Comparison between IHTMSr and IHTMS}\label{table:Num-res-all-comparison-r-IHTMS}
\begin{tabular}{c|c|c|c|c|c|c|c}\hline
\multicolumn{2}{c|}{Prob} & \multicolumn{3}{c|}{IHTMSr} &
\multicolumn{3}{c}{IHTMS} \\\hline

r &  FR  & NS & time & rel.err. & NS & time & rel.err.
\\\hline\hline

1    & 0.17      &10     & 2.59      & 1.67e-05         &10     &
3.98      & 1.77e-05 \\\hline

2    & 0.33      &10     & 4.98      & 2.11e-05         &10     &
6.95      & 2.04e-05 \\\hline

3    & 0.49      &10     & 9.95      & 2.27e-05         &10     &
12.65     & 2.30e-05 \\\hline

4    & 0.64      &10     & 22.72     & 3.05e-05         &10     &
27.12     & 2.86e-05 \\\hline

5    & 0.80      &10     & 84.13     & 3.95e-05         &10     &
94.13     & 4.10e-05 \\\hline

\end{tabular}
\end{center}
\end{table}

\begin{table}[ht]
\begin{center}\caption{Comparison between FPCAr and FPCA}\label{table:Num-res-all-comparison-r-FPCA}
\begin{tabular}{c|c|c|c|c|c|c|c}\hline
\multicolumn{2}{c|}{Prob} & \multicolumn{3}{c|}{FPCAr} &
\multicolumn{3}{c}{FPCA} \\\hline

r &  FR  & NS & time & rel.err. & NS & time & rel.err.
\\\hline\hline
1    & 0.17      &10     & 4.63      & 9.00e-06         &10     &
4.66      & 8.88e-06 \\\hline

2    & 0.33      &10     & 6.06      & 1.51e-05         &10     &
6.15      & 1.55e-05 \\\hline

3    & 0.49      &10     & 10.64     & 2.35e-05         &10     &
11.50     & 2.24e-05 \\\hline

4    & 0.64      &10     & 23.29     & 2.93e-05         &10     &
25.66     & 2.88e-05 \\\hline

5    & 0.80      &10     & 79.46     & 3.94e-05         &10     &
83.91     & 3.87e-05 \\\hline

\end{tabular}
\end{center}
\end{table}

From these tables we see that by using
our heuristics for determining the rank $r$ at every iteration, algorithms IHT, IHTMS and FPCA perform similarly to algorithms IHTr, IHTMSr and FPCAr which make use of knowledge of the true rank $r$. Specifically, algorithms IHT, IHTMS and FPCA are capable of recovering low-rank matrices very well even when we do not know their  rank.

\begin{table}[ht]
\begin{center}\caption{Comparison when the given rank is different from the true rank of 3}\label{table:wrong-rank}
\begin{tabular}{c|c|c|c}\hline

Given rank & NS & time & rel.err. \\\hline\hline

\multicolumn{4}{c}{IHTr} \\\hline

1 & 0 & --- & --- \\\hline

2 & 0 & --- & --- \\\hline

3 & 10  & 10.04 & 2.38e-05 \\\hline

4 & 10  & 21.42  & 3.42e-05   \\\hline

5 & 10  & 63.53 & 5.51e-05 \\\hline

6 & 4  & 109.00  & 4.44e-04  \\\hline\hline

IHT & 10 & 13.27  & 2.32e-05 \\\hline\hline

\multicolumn{4}{c}{IHTMSr} \\\hline

1 & 0 & --- & --- \\\hline

2 & 0 & --- & --- \\\hline

3 & 10  & 9.95 & 2.27e-05 \\\hline

4 &  10  & 22.53 & 3.40e-05   \\\hline

5 &  10 & 67.89 & 5.93e-05 \\\hline

6 &  1 & 116.62  & 6.04e-04  \\\hline\hline

IHTMS & 10 & 12.65 & 2.30e-05 \\\hline\hline

\multicolumn{4}{c}{FPCAr} \\\hline

1 & 0 & --- & --- \\\hline

2 & 0 & ---   & --- \\\hline

3 & 10  & 10.64 & 2.35e-05 \\\hline

4 & 10  & 21.26 & 3.46e-05 \\\hline

5 &  10  & 63.67 & 5.99e-05 \\\hline

6 &  3 & 108.02 & 4.04e-04 \\\hline\hline

FPCA & 10 & 11.50 & 2.24e-05 \\\hline\hline
\end{tabular}
\end{center}
\end{table}

Choosing $r$ is crucial in algorithms IHTr, IHTMSr and FPCAr as it
is in greedy algorithms for matrix rank minimization and compressed
sensing. In Table \ref{table:wrong-rank} we present results on how
the choice of $r$ affects the performance of algorithms IHTr, IHTMSr
and FPCAr when the true rank of the matrix is not known. In Table
\ref{table:wrong-rank}, the true rank is 3 and the results for
choices of the rank from 1 to 6 are presented. The rows labeled IHT,
IHTMS and FPCA present the results for these algorithms which use
the heuristics in Section \ref{sec:practical-issues} to determine
the rank $r$. From Table \ref{table:wrong-rank} we see that if
we specify a rank that is smaller than the true rank, then all of
the algorithms IHTr, IHTMSr and FPCAr are unable to successfully
recover the matrices (i.e., the relative error is greater than
1e-3). Specifically, since for the problems tested the true rank of
the matrix was 3, the algorithms failed when $r$ was chosen to be
either 1 or 2. If the chosen rank is slightly greater than the true
rank (i.e., the rank was chosen to be 4 or 5), all the three
algorithms IHTr, IHTMSr and FPCAr still worked. However, the
relative errors and times were much worse than those produced by the
heuristics based solvers IHT, IHTMS and FPCA. When the chosen rank
was too large (i.e., was chosen to be 6), IHTr, IHTMSr and FPCAr
were only able to recover the matrices in 4, 1 and 3 out of 10
problems, respectively. However, IHT, IHTMS and FPCA always recovered the matrices.

\subsection{A Video Compression Problem} We tested the performance of our algorithms on a video compression problem. By stacking each frame of the video as a column of a large matrix, we get a matrix $M$ whose $j$-th column corresponds to the $j$-th frame of the video. Due to the correlation between consecutive frames of the video matrix, $M$ is expected to be of low rank. Hence we should be able to recover the video by only taking a limited number of measurements. The video used in our experiment was downloaded from the website http://media.xiph.org/video/derf. The original colored video consisted of 300 frames where each frame was an image stored in an RGB format, as a $144\times 176\times 3$ array. Since this video data was too large for our use, we preprocessed it in the following way. We first converted each frame from an RGB format into a grayscale image, so each frame was a $144\times 176$ matrix. We then used only the portion of each frame corresponding to a $39\times 47$ submatrix of pixels in the center of each frame, and took only the first 20 frames. Consequently, the matrix $M$ had $m=1833$ rows and $n=20$ columns. We then created a Gaussian sampling matrix $A\in\br^{p\times (mn)}$ as in Section \ref{sec:numerical-results-random} with $p=1833*20*0.4=14664$ rows (i.e., we used sampling ratio $SR=0.4$) and computed $b=A\vvec(M)\in\br^p$. This $14664\times 36660$ matrix $A$ was close to the size limit of what could be created by calling the MATLAB function $A=randn(p,mn)$ on our computer. Although the matrix $M$ was expected to be of low rank, it was only approximately of low rank. Therefore, besides comparing the recovered matrices with the original matrix $M$, we also compared them with the best rank-$5$ approximation of $M$. Since the relative error of the best rank-$5$ approximation of $M$ was $2.33e-2$, we cannot expect to get a more accurate solution. Therefore, we set $xtol$ equal to $0.002$ for this problem. The results of our numerical tests are reported in Table \ref{table:Num-res-real}. The ranks reported in the table are the ranks of the recovered matrices. The reported relative errors and CPU times are averages over 5 runs. We do not report any results for SDPT3, because the problem is far too large to be solved by an SDP solver.
From Table \ref{table:Num-res-real} we see that our algorithms were able to recover the matrix $M$ very well, achieving relative errors that were of the same order as that obtained by the best rank-$5$ approximation.

\begin{table}[ht]
\begin{center}\caption{Results on recovery of compressed video}\label{table:Num-res-real}
\begin{tabular}{c|c|c|c}\hline
Solvers & rank & rel.err. & time  \\\hline

IHTr   & 5 & 6.87e-2 & 645 \\\hline

IHT    & 5 & 9.76e-2 & 949 \\\hline

IHTMSr & 5 & 6.72e-2 & 688 \\\hline

IHTMS  & 5 & 9.69e-2 & 804 \\\hline

FPCAr  & 5 & 5.10e-2 & 514 \\\hline

FPCA   & 5 & 5.17e-2 & 1296 \\\hline

\end{tabular}
\end{center}
\end{table}

\begin{figure}
\begin{center}
\includegraphics[scale=0.65]{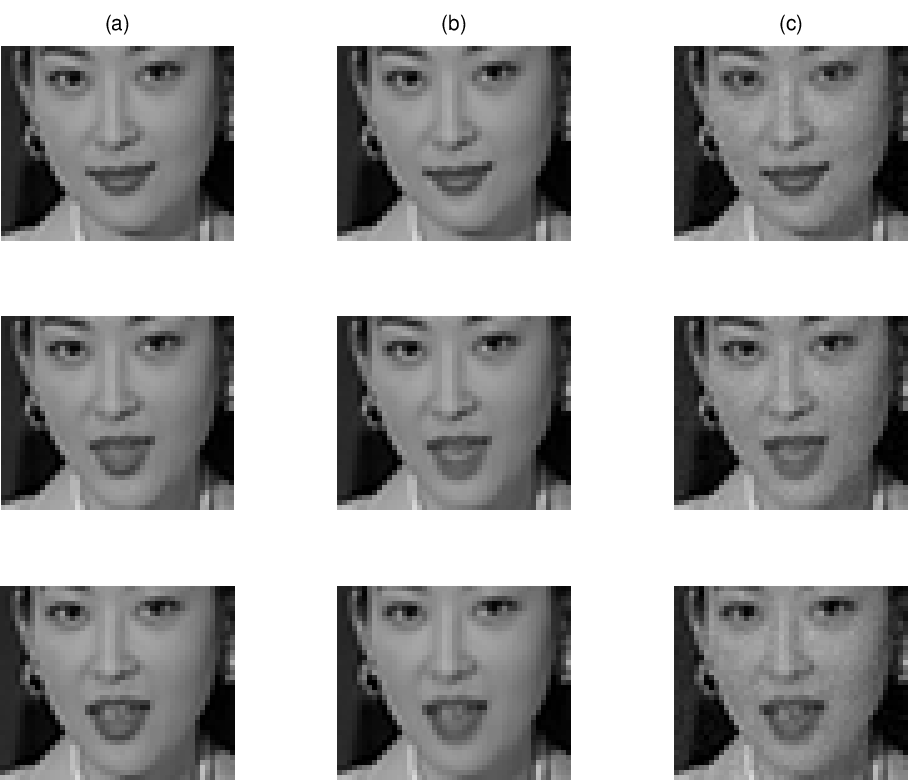}
\caption{Comparison of frames 4, 12 and 18 of (a) the original video, (b) the best rank-5 approximation and (c) the matrix recovered by FPCA} \label{fig:video}
\end{center}
\end{figure}

In Figure \ref{fig:video}, the three images in the first column correspond to three particular frames in the original video. The images in the second column correspond to these frames in the rank-$5$ approximation matrix of the video. The images in the third column correspond to these frames in the matrix recovered by FPCA. The other five solvers recovered images that were very similar visually to FPCA so we do not show them here. From Figure \ref{fig:video} we see that FPCA recovers the video very well by taking only 40\% as many measurements as there are pixels in the video.

\section*{Acknowledgement}  We would like to thank two anonymous referees for insightful comments that greatly improved the presentation of the paper. We would also like to thank Dr. Thomas Blumensath for pointing out an error in an earlier version of this paper.

\renewcommand{\theequation}{A-\arabic{equation}}
\setcounter{equation}{0}  
\renewcommand{\thetheorem}{A-\arabic{theorem}}
\setcounter{theorem}{0}
\section*{Appendix}  

Here we give proofs of Propositions \ref{prop:RIP:PA*,PA*AP},
\ref{prop:RIP:PPsi1A*APPsi2} and \ref{prop:RIP:error bound of AX}.

{\bf Proof of Proposition \ref{prop:RIP:PA*,PA*AP}}. \proof We
prove \eqref{prop:RIP:PA*,PA*AP-eq1} first. Since for any
$X\in\br^\mtn,$
$\rank(P_\Psi X)\leq r$, we have
\begin{align*}\begin{split}\left|\langle X,P_\Psi\A^* b \rangle\right| & = |\langle\A P_\Psi X,b\rangle|\\
                                             & \leq \|\A P_\Psi X\|_2\|b\|_2\\
                                             & \leq \sqrt{1+\delta_r(\A)}\|P_\Psi X\|_F\|b\|_2\\
                                             & \leq \sqrt{1+\delta_r(\A)}\|X\|_F\|b\|_2. \end{split}\end{align*}
Thus \begin{align*}\|P_\Psi\A^* b\|_F =
\max_{\|X\|_F=1}|\langle X,P_\Psi\A^* b \rangle| \leq \sqrt{1+\delta_r(\A)}\|b\|_2.\end{align*}

To prove \eqref{prop:RIP:PA*,PA*AP-eq2}, note that by the RIP,
\begin{align*}(1-\delta_r(\A))\|P_\Psi X\|_F^2\leq \|\A P_\Psi X\|_F^2 \leq (1+\delta_r(\A))\|P_\Psi X\|_F^2,\end{align*} which means the eigenvalues of $P_\Psi\A^*A
P_\Psi$ restricted to $\myspan(\Psi)$ are in the interval
$[1-\delta_r(\A),1+\delta_r(\A)]$. Thus
\eqref{prop:RIP:PA*,PA*AP-eq2} holds.
 \eproof

{\bf Proof of Proposition \ref{prop:RIP:PPsi1A*APPsi2}}.
First, we prove
\begin{align}\label{proof-Prop-PA*A(I-P)}
\left|\langle\A(I-P_\Psi)X, \A P_\Psi Y \rangle\right| \leq \delta_{r}(\A)\|(I-P_\Psi)
X\|_F\|P_\Psi Y\|_F, \forall Y\in\mathbb{R}^\mtn,
X\in\myspan(\Psi').
\end{align}
\eqref{proof-Prop-PA*A(I-P)} holds obviously if $(I-P_\Psi) X =0$ or
$P_\Psi Y = 0$. Thus we can assume $(I-P_\Psi) X \neq 0$ and $P_\Psi
Y \neq 0.$ Define $\hat{X}=\frac{(I-P_\Psi) X}{\|(I-P_\Psi) X\|_F}$
and $\hat{Y} = \frac{P_\Psi Y}{\|P_\Psi Y\|_F}$; then we have
$\left\|\hat{X}\right\|_F = 1$, $\left\|\hat{Y}\right\|_F=1$ and $\langle\hat{X}, \hat{Y}\rangle=0.$
Since $\hat{X}\in\myspan(\Psi\cup\Psi')$ and
$\hat{Y}\in\myspan(\Psi)$, we have $\rank\left(\hat{X}+\hat{Y}\right)\leq r$
and $\rank\left(\hat{X}-\hat{Y}\right)\leq r$. Hence by RIP,
\begin{align*}\begin{split}
2(1-\delta_{r}(\A))= (1-\delta_{r}(\A))\left\|\hat{X}+\hat{Y}\right\|_F^2&\leq
\left\|\A \hat{X} + \A \hat{Y} \right\|_2^2 \\ & \leq
(1+\delta_{r}(\A))\left\|\hat{X}+\hat{Y}\right\|_F^2=2(1+\delta_{r}(\A)).
\end{split}
\end{align*}
and
\begin{align*}\begin{split}
2(1-\delta_{r}(\A))= (1-\delta_{r}(\A))\left\|\hat{X}-\hat{Y}\right\|_F^2 &
\leq \left\|\A \hat{X} - \A \hat{Y} \right\|_2^2 \\ & \leq
(1+\delta_{r}(\A))\left\|\hat{X}-\hat{Y}\right\|_F^2=2(1+\delta_{r}(\A)).
\end{split}
\end{align*}
Therefore we have
\begin{align*}\langle\A \hat{X},\A \hat{Y}\rangle = \frac{\left\|\A \hat{X}+\A \hat{Y}\right\|_2^2-\left\|\A \hat{X}- \A \hat{Y} \right\|_2^2}{4}\leq \delta_{r}(\A) \end{align*}
and
\begin{align*}-\langle\A \hat{X},\A \hat{Y}\rangle = \frac{\left\|\A \hat{X}-\A \hat{Y}\right\|_2^2-\left\|\A \hat{X}+ \A \hat{Y} \right\|_2^2}{4}\leq \delta_{r}(\A). \end{align*}
Thus, $|\langle\A \hat{X},\A \hat{Y}\rangle|\leq\delta_{r}(\A)$ and
\eqref{proof-Prop-PA*A(I-P)} holds.

Finally we have, for any $X\in \myspan(\Psi')$,
\begin{align*}\begin{split}
\left\|P_\Psi\A^*\A(I-P_\Psi)X\right\|_F & = \max_{\|Y\|_F=1}|\langle P_\Psi
\A^*\A(I-P_\Psi)X, Y\rangle| \\ & = \max_{\|Y\|_F=1}|\langle\A(I-P_\Psi)X, \A
P_\Psi Y\rangle| \\ & \leq \delta_{r}(\A)\left\|(I-P_\Psi)X\right\|_F, \end{split}
\end{align*} i.e., \eqref{prop:RIP:PPsi1A*APPsi2-result} holds,
which completes the proof.
 \eproof

{\bf Proof of Proposition \ref{prop:RIP:error bound of AX}}. This
proof essentially follows that given by Needell and Tropp in
\cite{Needell-Tropp-2008}. \proof Let $B^s:=\{X\in\br^\mtn:
\rank(X)=s,\|X\|_F\leq 1\}$ be the unit ball of rank-$s$ matrices in
$\br^\mtn.$ Define the convex hull of the unit norm matrices with
rank at most $r$ as:
\begin{align*}S:=\mbox{conv}\left\{\bigcup_{s\leq r}B^s\right\}\subset\br^\mtn.\end{align*}
By \eqref{prop:RIP:error bound of AX-hypothesis-eq}, we know that
the operator norm
\begin{align*}\|\A\|_{S\rightarrow 2}=\max_{X\in S}\|\A X\|_2\leq \sqrt{1+\delta_r(\A)}.\end{align*}
Define another convex set
\begin{align*}K:=\{X\in\br^\mtn:\|X\|_F+\frac{1}{\sqrt{r}}\|X\|_*\leq 1\}\subset\br^\mtn,\end{align*}
and consider the operator norm
\begin{align*}\|\A\|_{K\rightarrow 2} = \max_{X\in K}\|\A
X\|_2.\end{align*} The content of the proposition is the claim that
$K\subset S.$

Choose a matrix $X\in K$ with SVD $X=U\Diag(\sigma)V^\top$. Let $I_0$ index the $r$
largest components of $\sigma$, breaking ties lexicographically. Let
$I_1$ index the next largest $r$ components, and so forth. Note that
the final block $I_J$ may have fewer than $r$ components. We may
assume that $\sigma|_{I_j}$ is nonzero for each $j$. This partition
induces a decomposition
\begin{align*}X=U[\Diag(\sigma|_{I_0})+\sum_{j=1}^J \Diag(\sigma|_{I_j})]V^\top=\lambda_0 Y_0 + \sum_{j=1}^J\lambda_j Y_j,\end{align*}
where $\lambda_j=\|U\Diag(\sigma|_{I_j})V^\top\|_F$ and $Y_j =
\lambda_j^{-1}U\Diag(\sigma|_{I_j})V^\top$. By construction, each
matrix $Y_j$ belongs to $S$ because it's rank is at most $r$ and it
has unit Frobenius norm. We will prove that $\sum_{j}\lambda_j\leq
1,$ which implies that $X$ can be expressed as a convex combination
of matrices from the set $S$. So $X\in S$ and $K\subset S.$

Fix $j$ in the range $\{1,2,\ldots, J\}.$ It follows that
$\sigma|_{I_j}$ contains at most $r$ elements and
$\sigma|_{I_{j-1}}$ contains exactly $r$ elements. Therefore,
\begin{align*}\lambda_j=\|\sigma|_{I_j}\|_2\leq\sqrt{r}\|\sigma|_{I_j}\|_\infty\leq\sqrt{r}\cdot\frac{1}{r}\|\sigma|_{I_{j-1}}\|_1.\end{align*}
Summing these relations, we obtain,
\begin{align*}\sum_{j=1}^J\lambda_j\leq\frac{1}{\sqrt{r}}\|\sigma|_{I_{j-1}}\|_1\leq\frac{1}{\sqrt{r}}\|X\|_*.\end{align*} It is obvious that
$\lambda_0=\|\sigma|_{I_0}\|_2\leq \|X\|_F.$ We now conclude that
\begin{align*}\sum_{j=0}^J\lambda_j\leq \|X\|_F+\frac{1}{\sqrt{r}}\|X\|_*\leq 1\end{align*} because $X\in K.$
This implies that $X\in S$ and $K\subset S$, and thus completes the
proof. \eproof

\bibliographystyle{acm}
\bibliography{All}
\end{document}